%% file: main_ArXiv_v2.tex
\documentclass[preprint,12pt]{elsarticle}
\usepackage[left=1.0in,right=1.0in,top=1in,bottom=1in]{geometry}

\usepackage{subfiles}
\usepackage{float}
\usepackage{ulem}

\usepackage{graphics,graphicx}
\usepackage{setspace} 
\usepackage{wrapfig}
\usepackage[dvipsnames,x11names]{xcolor}
\usepackage{longtable}
\usepackage{caption}
\usepackage{subcaption}
\usepackage{amsmath}
\usepackage{amsfonts}
\usepackage{amssymb}
\usepackage{placeins}
\usepackage{rotating}
\usepackage{bm}
\usepackage[hyphens]{url}
\usepackage[hidelinks]{hyperref}
\usepackage{mathnotation}
\usepackage{lscape}
\usepackage{multirow}
\usepackage{booktabs}
\usepackage{stmaryrd}
\usepackage{tikz}

\newcommand{\RevOne}[1]{\textcolor{black}{#1}}
\newcommand{\RevTwo}[1]{\textcolor{black}{#1}}
\newcommand{\RevThree}[1]{\textcolor{black}{#1}}

\newcommand\includegraphicsifexists[2][width=\linewidth]{\IfFileExists{#2}{\includegraphics[#1]{#2}}{}}
\graphicspath{{Figures}}

\newcommand{\avg}[1]{\{\!\{#1\}\!\}}

\journal{Journal of Computational Physics}

\input{macros/FXG_Macros}

\begin{document}

\begin{frontmatter}
\title{GPU Performance of an Entropy-Stable Discontinuous Galerkin Euler Solver with Non-Conservative Terms}

\author[ucsc]{Henry Waterhouse}
\author[snl]{Maciej Waruszewski}
\author[nps]{Lucas C. Wilcox}
\author[vatech]{Timothy Warburton}
\author[ucsc,nps]{Francis X. Giraldo \corref{cor1} }
\cortext[cor1]{Corresponding author}

\address[ucsc]{Department of Applied Mathematics, University of California, Santa Cruz, CA}
\address[snl]{Sandia National Laboratories, Albuquerque, NM}
\address[nps]{Department of Applied Mathematics, Naval Postgraduate School, Monterey, CA}
\address[vatech]{Department of Mathematics, Virginia Tech, Blacksburg, VA}

\begin{abstract}
The entropy-stable discontinuous Galerkin method for the compressible Euler equations with buoyancy is implemented on graphics processing unit (GPU) hardware. We measure the performance of the solver on three-dimensional problems: the rising thermal bubble and the baroclinic instability in a channel. On NVIDIA A100 hardware, the solver achieves nearly 70\% of 64-bit floating-point peak performance for the most computationally expensive kernel (volume terms) and significantly reduces the computational overhead typically incurred by two-point entropy-stable fluxes in the volume terms. 
We also demonstrate strong and weak scaling performance of the solver and 
compare it with a highly optimized central processing unit (CPU) code, showing that
the GPU kernels are 
$10\times$ faster
and $13\times$ more energy efficient than the CPU code. We also show that the solver achieves the expected $2\times$ speedup when run at 32-bit floating-point precision. We discuss the 
modifications used to obtain the final GPU implementation
and measure the performance gain of each of the implementation 
strategies, including reductions in complex operations and memory traffic and improvements in load balancing.
We also extend symmetry-based flux savings to the non-symmetric gravity term, preserving nearly the full factor-of-two speedup achieved for the symmetric flux.
\end{abstract}

\begin{keyword}
atmospheric model \sep balance laws \sep baroclinic instability \sep buoyancy \sep compressible Euler equations \sep discontinuous Galerkin \sep entropy-stable method \sep flux-differencing \sep nonhydrostatic \sep numerical weather prediction \sep rising thermal bubble \sep Runge-Kutta method
\end{keyword}

\end{frontmatter}

\input{sections/Introduction}

\input{sections/GovernEq}

\input{sections/NumericalMethod}

\input{sections/gpu_optimizations}
\input{sections/performance_analysis}

\input{sections/simulation_results}

\input{sections/Conclusions}

\section*{Acknowledgments}  \label{sec:acknowledgement}
Waruszewski, Wilcox, and Giraldo were supported by the generosity of Eric and Wendy Schmidt by recommendation of the Schmidt Futures program. Waruszewski, Wilcox, and Giraldo were also supported by the National Science Foundation under grant AGS-1835881. Giraldo was also supported by the Office of Naval Research Marine Meteorology and Space Weather program under grant N0001419WX00721.
Sandia National Laboratories is a multimission laboratory managed and operated by National Technology \& Engineering Solutions of Sandia, LLC, a wholly owned subsidiary of Honeywell International Inc., for the U.S. Department of Energy’s National Nuclear Security Administration under contract DE-NA0003525.
This paper describes objective technical results and analysis. Any subjective views or opinions that might be expressed in the paper do not necessarily represent the views of the U.S. Department of Energy or the United States Government.
This work used Delta at the National Center for Supercomputing Applications through allocation MTH250049 from the Advanced Cyberinfrastructure Coordination Ecosystem: Services \& Support (ACCESS) program, which is supported by U.S. National Science Foundation grants \#2138259, \#2138286, \#2138307, \#2137603, and \#2138296.


\bibliographystyle{siam}
\bibliography{bibliography/henry_references, bibliography/Giraldo_BibFile_April_14_2026}

\end{document}

%% file: macros/FXG_Macros.tex
\usepackage[english]{babel} 

\usepackage{epsfig}

\usepackage{algorithmicx}
\usepackage{algpseudocode}

\usepackage{stmaryrd}
\usepackage{amsmath} 
\usepackage{color}
\usepackage{amssymb} 
\usepackage{amsfonts} 
\usepackage{mathrsfs}

\usepackage[T1]{fontenc}
\usepackage[utf8]{inputenc}
\usepackage[font=small,labelfont=bf,tableposition=top]{caption}
\usepackage{multicol}

\usepackage{float}
\usepackage{graphicx}
\input{epsf}

\newcommand{\comment}[1] {}

\newcommand{\diff}[2] {\frac{\partial #1}{\partial #2}}

\newcommand{\vc}[1]{\mbox{\boldmath$#1$\unboldmath}}

\newcommand{\be}{\begin{equation}}
\newcommand{\ee}{\end{equation}}
\newcommand{\bea}{\begin{eqnarray*}}
\newcommand{\eea}{\end{eqnarray*}}

\newcommand{\floor}[1]{\left\lfloor #1 \right\rfloor}

%% file: sections/Introduction.tex

\section{Introduction} 
\label{sec:introduction}

The construction of accurate, stable, and efficient dynamical cores in atmospheric modeling remains a topic of interest despite the volume of literature on this subject. Most (perhaps all) existing atmospheric models require some form of regularization to maintain stability. Typically, this regularization has come in the form of hyper-diffusion (see, e.g., \cite{Giraldo1999, Kim2008}). Although this idea has been improved to make it more adaptive (see, e.g., \cite{Guba2014a, Giraldo2024}), some tuning is still required when the model resolution is modified. Discontinuous Galerkin (DG) methods offer many advantages in geophysical fluid dynamics modeling 
because of their 
(1) high-order accuracy (see, e.g., \cite{Hesthaven2008, Kopriva2009, Giraldo2020}), (2) stability due to the use of Riemann solvers (essentially augmenting the numerics with analytic information from the method of characteristics, see, e.g., \cite{toro_riemann_2009}), and (3) the local element-wise construction, which offers unmatched parallelization (i.e., efficiency) \cite{Kelly2012, Abdi2017a, Abdi2017b}.  The main limitations of these methods have been (i) stability at high order and (ii) time-step restrictions.  The stability issue has been resolved thanks to entropy-stable methods (see, e.g., \cite{Tadmor1984, Tadmor2006, Castro2013, Fisher2013, Kopriva2014, Gassner2016, Kopriva2016a, Parsani2016, Chen2017, Friedrich2018, Chan2019, Chan2019a,Renac2019,Waruszewski2022,Souza2023}) which seek to construct the correct numerical flux that either conserves entropy (entropy-conservative methods) or dissipates it (entropy-dissipative or entropy-stable methods); this is achieved by deriving the numerical flux for the associated entropy equation that achieves these goals. The time-step issue has been ameliorated, to some degree, by implicit methods (see, e.g., \cite{Reddy2023, Giraldo2024, Welter2026}). 
These two approaches address the stability and time-step restrictions,
leaving the cost of entropy-stable methods to be addressed. 
Entropy-stable discontinuous Galerkin (ESDG) methods are expensive to compute due to the use of two-point fluxes (see \cite{Fisher2013}). Therefore, it is the main objective of this paper to show how to construct efficient ESDG methods on modern high-performance computing architectures, i.e., graphics processing units (GPUs); Ranocha et al.\ \cite{Ranocha2023} tackled many of these issues on the central processing unit (CPU).  We aim to show how each optimization we apply improves the efficiency, thereby developing a roadmap for constructing efficient GPU code for various applications, although here we focus on the Euler equations with buoyancy. The challenge of including buoyancy in the Euler equations is that the resulting system of nonlinear partial differential equations no longer forms a  system of conservation laws but rather a system of balance laws with non-conservative terms that must be handled carefully; to this end, we follow the ESDG method presented in \cite{Waruszewski2022}. To our knowledge, the only published works on the GPU performance of ESDG have focused exclusively on the shallow water equations \cite{Wintermeyer2018, Wu2021}; our manuscript represents the first work on ESDG for the Euler equations that focuses on GPU optimization---our  aim is to demonstrate 
the effects on performance of various optimizations, such as reductions in complexity and memory usage.

It is important to note that, in the last two decades, 
a substantial body of work has examined GPU implementations
for the Euler equations, e.g.,  \cite{Elsen2008,Kuo2011,Lefebvre2012,Siebenborn2012,Siebenborn2013,Fuhry2014,Witherden2014,Abdi2017a, Abdi2017b,Zhang2018,Giuliani2019,Kirby2020,Romero2020,Brodtkorb2022,Wang2022,Cernetic2023,Cernetic2024,Kurz2025,Witherden2025,Sporykhin2026,Welter2026}. A smaller collection focuses on the GPU implementation of discontinuous Galerkin methods for the shallow water equations \cite{Wintermeyer2018,Wu2021} and the Euler equations \cite{Siebenborn2012,Siebenborn2013,Fuhry2014,Abdi2017a, Abdi2017b,Kirby2020,Wang2022,Cernetic2023,Cernetic2024,Kurz2025,Welter2026}. From this list, only two papers focus on the GPU implementation of ESDG methods for the shallow water equations (SWE), and none focuses on the Euler equations.\footnote{\RevThree{Comparing results between SWE and Euler is difficult because the 
ESDG kernels have much lower complexity and cost for SWE than for Euler, and the SWE results are two-dimensional, whereas the results presented here are three-dimensional.}}
This literature review shows that although a rather large volume of work has appeared on ESDG methods, 
very few studies have addressed GPU implementations of these methods.  Our goal is to fill this important gap in the literature because, without improving performance, ESDG methods (as good as they are) run the risk of falling out of favor \RevThree{due to} their computational cost.

The remainder of the paper is organized as follows: Sec.~\ref{sec:governing_equations} describes the form of the Euler equations used in this study, which includes the effect of buoyancy, and Sec.~\ref{sec:numerical_method} presents the ESDG discretization. 
Section~\ref{sec:performance} presents the performance analysis: Sec.~\ref{sec:optimization} describes the optimizations, and Secs.~\ref{sec:single-gpu} and~\ref{sec:multi-gpu} present the single- and multi-GPU results, respectively.
We test our algorithm on three-dimensional test cases in Sec.~\ref{sec:simulations} and conclude the paper in Sec.~\ref{sec:conclusions}.

%% file: sections/GovernEq.tex
\section{Governing Equations} 
\label{sec:governing_equations}

This study considers flow governed by the time-dependent compressible Euler equations. Let \RevThree{$\Omega \subset \mathbb{R}^3$} be a fixed domain in 3D with boundary $\Gamma$. Let us denote the fluid density, velocity vector, and total energy by $\rho, \vc{u}=u_i$ $(i=1,2,3)$, and $e$, where $e=c_v T + \frac{1}{2}\vc{u}\cdot \vc{u} + \phi$. 
Here, \(c_v\) denotes the specific heat at constant volume, \(T\) denotes temperature, and \(\phi=gz\) denotes the geopotential, with \(g\) denoting the gravitational constant.
The Euler equations are defined in conservation form as follows:
\subeqs{Eq:euler}{
\label{Eq:euler}
\eq{mass_eq}{ \diff{\rho}{t} + \nabla \cdot \vc{U} = 0 }
\eq{mom_eq}{\diff{\vc{U}}{t} + \nabla \cdot \left[ \frac{1}{\rho} \vc{U} \otimes \vc{U} + P \vc{I}_3  \right] + \rho \nabla \phi + \vc{\omega} \times \vc{U} = 0 } 
\eq{energy_eq}{ \diff{E}{t} + \nabla \cdot \left[ \left( \frac{E+P}{\rho}\right) \vc{U} \right] =  0 }
}
where $\vc{U}=\rho \vc{u}$ is the momentum, $E=\rho e$ is the total energy density (which includes kinetic, internal, and potential energies), $\omega$ is the rotation rate contributing to the Coriolis effect, P is the pressure defined as
\expression{
P=(\gamma -1)\left(E - \frac{1}{2 \rho} \vc{U}\cdot\vc{U} - \rho \phi \right),
}
where $\gamma=\frac{c_p}{c_v}$ is the specific heat ratio, and $\vc{I}_3$ is the rank-3 identity matrix. 

A few remarks on Eq.~(\ref{Eq:euler}) are warranted. Ignoring the geopotential (buoyancy) and Coriolis terms in momentum results in a conservation law. 
The buoyancy and Coriolis terms are what distinguish the standard Euler equations in computational fluid dynamics from those required in geophysical fluid dynamics and result in a balance law. It is the non-conservative buoyancy term in particular that causes us to construct specialized entropy-stable discretizations as outlined in the next section.

%% file: sections/NumericalMethod.tex
\section{Numerical Method} 
\label{sec:numerical_method}

The ESDG discretization this study considers
is detailed in \cite{Waruszewski2022}. Here we provide a brief description of
the method with a focus on the operations the computational kernels implement.

The domain is divided into a finite number of non-overlapping
elements which make up the computational domain. Each element is
mapped to the hexahedral reference element containing a tensor-product grid of
$(N+1)^3$ Legendre-Gauss-Lobatto (LGL) quadrature nodes, where $N$ is the
polynomial order of approximation. The discrete state vector stores the
approximate solution at these nodes and all integrals are evaluated using the
associated quadrature weights (to increase computational efficiency by avoiding
projections between interpolation and integration nodes). For each element, the
semi-discrete form of the ESDG discretization
is given as
\begin{equation}\label{nm:semi}
\frac{d\mathrm{q}}{dt} =
\underbrace{%
\mathrm{h} 
- \sum_{j,k=1}^3
  \left(\mathrm{G}_{jk} \left(\mathrm{D}_j \circ \mathrm{F}_j \right)
- \left(\mathrm{F}_j \circ \mathrm{D}^T_j\right) \mathrm{G}_{jk} \right)
\mathrm{\mathbf{1}}}_{\text{volume term}}
-
\underbrace{%
\sum_{j,k=1}^3
\mathrm{M}^{-1}
\mathrm{M}^s
\mathrm{G}^s_{jk}
\mathrm{F}^*_j}_{\text{surface term}},
\end{equation}
where \(\mathrm{q}\) is a vector of the state approximated at the nodes on
the element;
\(\mathrm{h}\) is a vector containing a source evaluated at the nodes;
\(\mathrm{G}_{jk}\) is 
the diagonal matrix of geometric factors
containing the metric
derivatives of the $k$-th direction reference coordinate with respect to the
$j$-th direction spatial coordinate;
\(\mathrm{D}_k\) is the $k$-th reference direction nodal differentiation
matrix;
\(\mathrm{F}_j\) is the $j$-th spatial direction two-point flux matrix that ensures entropy conservation and \(\mathrm{F}^*_j\) is the corresponding entropy-dissipative flux defined as the entropy-conservative flux with additional dissipation included;
\(\mathrm{\mathbf{1}}\) is a vector of ones;
\(\mathrm{M}\) is the Jacobian-weighted LGL-quadrature-based diagonal element
mass matrix;
\(\mathrm{G}^s_{jk}\) is 
the diagonal matrix of geometric factors
for the element
surface nodes; and
\(\mathrm{M}^s\) is the Jacobian-weighted LGL-quadrature-based diagonal surface
mass matrix.
Volume and surface terms are applied in a matrix-free way, exploiting the
tensor-product structure of the differentiation matrices. For example, only an
$N_q \times N_q$ matrix will need to be loaded to apply each derivative, where
$N_q=N+1$ is the number of quadrature points in each dimension of the reference
element. The elements are connected together through the calculation of the
surface flux vector, which requires state values from the current element (the
minus side of the surface) and the neighboring elements (the plus side of the
surface). Entropy stability is ensured by evaluating the flux using a specific
two-point flux formula which allows one to prove discrete entropy
stability~\cite{Waruszewski2022}.
For the volume kernels, we use the entropy-conserving two-point flux described in \cite{Waruszewski2022} (Eq.\ 76); for the surface flux we use the same entropy-conserving flux in addition to the matrix dissipation outlined in \cite{Waruszewski2022} (Eq.\ A.1). The entropy-conserving flux is computationally intensive because its construction requires mean values, jump values, logarithms, and reciprocals. The matrix dissipation is somewhat computationally \RevThree{intensive} because it requires diagonalization of the Jacobian matrix, akin to a Roe Riemann solver \cite{Roe1986}.
As detailed in Sec.~\ref{sec:performance}, this is done by
exploiting symmetry and the tensor-product structure to minimize the number of
computations.

Although we previously relied on the 5-stage fourth-order low-storage Runge--Kutta scheme~\cite{CarpenterKennedy1994}, we present results using a  fourteen-stage fourth-order low-storage Runge--Kutta
scheme~\cite{Niegemann2012} due to the decrease in time-to-solution that it has provided us. The elements of the two-point flux matrices
\(\mathrm{F}_j\) and vectors \(\mathrm{f}_j\) depend on the state and need to
be recomputed at each stage. This accounts for a significant portion of the
floating-point operations performed in the volume and surface kernels.

%% file: sections/gpu_optimizations.tex
\section{Performance Analysis}
\label{sec:performance}

\subsection{\RevOne{Algorithmic and GPU implementation optimizations}}
\label{sec:optimization}

DG methods operate on a set of independent elements that only communicate through surface fluxes at the boundaries.
Most of the work is done inside an element, in what is usually referred to as the volume term, which is well suited to GPU optimizations.
The GPU thread hierarchy is organized into thread blocks.
Threads within a block can cooperate, but each thread block is scheduled independently. The first step in GPU optimization is to determine an efficient mapping of computational work to GPU threads.
DG methods provide an obvious distribution methodology: each element is represented by a thread block, with each thread within that block representing a single degree of freedom (DOF) within the computational element.

\begin{figure}[htbp]
    \centering
    \includegraphics[width=0.8\textwidth]{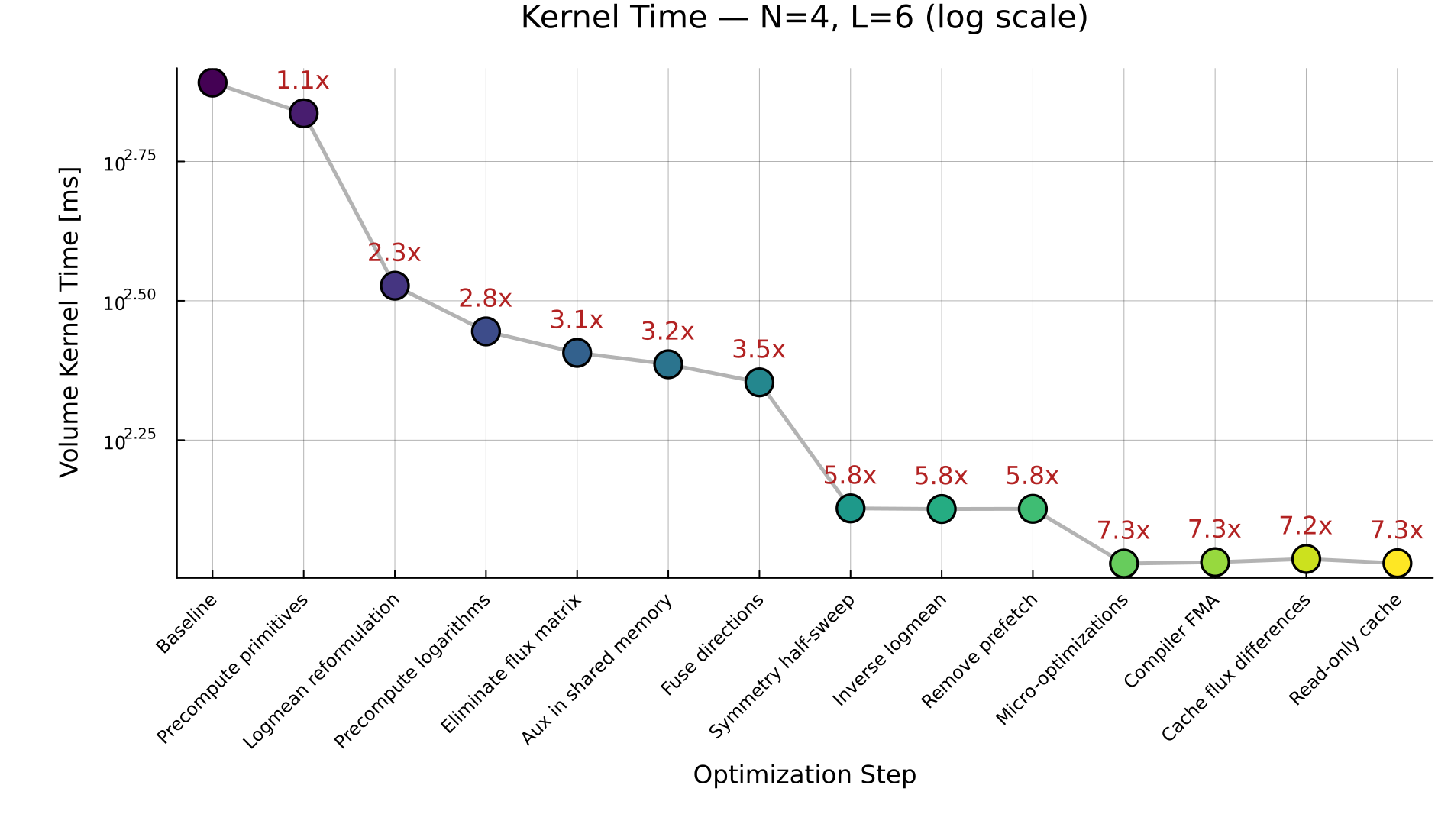}
    \caption{Total volume kernel GPU time per timestep for each cumulative optimization step. Measured on an NVIDIA A100 SXM4 40GB with polynomial order $N=4$, refinement level $L=6$, 64-bit precision, using NVIDIA Nsight Compute. Although we used the rising thermal bubble described in Sec.~\ref{sec:rising-bubble-case} as the underlying simulation, the results should be test-case agnostic.}
    \label{fig:kernel-timings}
\end{figure}

Along with this distribution of work, there are many optimizations to be made to the kernel itself.
We build on the efficient two-point flux formulation of \cite{Ranocha2023}, which benchmarked serial central processing unit (CPU) performance, and extend it to the GPU. 
\RevTwo{Our code is written in the Julia programming language \cite{bezanson2017julia} and uses the package CUDA.jl \cite{besard2018juliagpu} to implement the GPU kernels discussed below.}
Figure~\ref{fig:kernel-timings} shows the cumulative effect of each optimization step, culminating in a total speedup of $7.3\times$ over the baseline GPU kernel.\footnote{\RevThree{The kernel implementation corresponding to each optimization step in Fig.~\ref{fig:kernel-timings} is available at \url{https://github.com/HorribleSanity/ESDGPerf}}.}
In Fig.~\ref{fig:kernel-timings}, $N$ denotes the polynomial degree of the basis function for each direction and the number of elements in the mesh is defined by $N_e=K \times 2^{3L}$ where $K=K_x \times K_y \times K_z$ and $K_j$ denotes the number of base elements along direction $j$ and $L$ denotes the refinement level of the mesh.
Unless specifically stated, test cases are run with only one base element ($K = 1$).
We discuss the most significant optimizations below.

The two-point flux used in flux-differencing methods (e.g., entropy-stable and kinetic-energy-preserving methods) is evaluated between every pair of quadrature nodes along each coordinate direction, requiring $O(N_q)$ evaluations per DOF per direction. 
The high number of evaluations of this flux function makes it crucial to avoid as many expensive operations within the function as possible.
Two classically expensive operations in high-performance computing are divisions and logarithms, both of which feature prominently in entropy-stable fluxes.
Logarithms are particularly costly in 64-bit precision on GPUs, where transcendental functions lack dedicated hardware support \cite{nvidia_cuda_nodate}.
A major source of division operations is the computation of primitive variables
(e.g., $u = \rho u / \rho$).
We limit this bottleneck by precomputing these once per thread instead of as needed in the flux.
The logarithms appear in the form of logarithmic means.
We adopt the efficient logarithmic mean formulation of \cite{Ranocha2023}, which reduces divisions with Horner's rule.
We further rewrite the logarithmic mean so that it depends on $\log(a)$ and $\log(b)$ individually rather than on $\log(\frac{a}{b})$, allowing the logarithms to be computed once per thread and reused across all evaluations.
These changes make each evaluation of the two-point flux cheaper, but we can also lower the total number of two-point fluxes by using symmetry.

\RevOne{We refer to exploiting symmetry to evaluate the entropy-stable flux once per unordered pair of quadrature nodes along each coordinate direction, rather than once for each ordering of the pair, as a ``half-sweep''.}
\RevOne{The half-sweep optimization for} the standard Euler equations, where the entropy-stable flux is built from symmetric averages, is described in \cite{Ranocha2023}.
Here, we encounter two additional complications. The first one is that the gravity term results in  non-symmetric terms added to the entropy-stable flux. The second is that the standard CPU way of computing only the necessary two-point fluxes leads to load imbalance on GPUs. The first problem comes from the momentum flux. We show this by first defining:
\begin{align}
    \avg{a} = \frac{a^+ + a^-}{2}, \,\,\,\, \avg{a}_{\mathrm{log}} = \frac{a^+ - a^-}{\log(a^+) -\log(a^-)}, \,\,\,\, \llbracket a \rrbracket = a^+ - a^-
\end{align}
and then stating the form of the momentum flux:
\begin{align}
    \avg{\rho}_{\mathrm{log}} \avg{u_k} \avg{u_m} + \delta_{mk}\left(\tfrac{\avg{\rho}}{2\avg{b}} + \tfrac{1}{2}\hat{\rho} \llbracket \phi \rrbracket \right),
\end{align}
where
\begin{align}
    \hat{\rho} = \frac{\avg{b} \avg{\rho}_{\mathrm{log}}}{b^-}
\end{align}
contains the one-sided expression $b^- = \rho^-/(2p^-)$.
Under an $i \leftrightarrow j$ swap, the jump $\llbracket\phi\rrbracket$ flips sign and $b^-$ becomes $b^+$, so the gravity term transforms as $\tfrac{1}{2}\hat{\rho}\llbracket\phi\rrbracket \to -\tfrac{1}{2}\hat{\rho}\llbracket\phi\rrbracket \cdot \frac{b^-}{b^+}$, which is not a simple sign flip when $b^- \neq b^+$. We handle this by decomposing the momentum flux into a symmetric part $S = \avg{\rho}_{\mathrm{log}} \avg{u_k}\avg{u_m} + \delta_{mk} \tfrac{\avg{\rho}}{2\avg{b}}$ and the gravity contribution $G = \tfrac{1}{2}\delta_{mk}\hat{\rho}\llbracket\phi\rrbracket$.
For each quadrature node pair, the flux is computed once and the partner node's gravity contribution is obtained by a single scaling: $G \to -G \cdot \tfrac{b^-}{b^+}$.
The density and energy components remain fully symmetric and require only a sign flip.
This preserves nearly all of the factor-of-two savings from the half-sweep while correctly handling the gravitational coupling.
The second challenge of exploiting symmetry comes from GPU threading.
A naive assignment of work where thread $i$ computes two-point fluxes $F_{i,j}$ for $ j > i$ leads to load imbalance.
We can handle this by making thread $i$ compute $F_{i, 1 + (i + l - 2) \bmod N_q}$ for $l = 2, \dots , N_l$ where
\begin{align}
    N_l = \begin{cases}
    \floor{\frac{N_q}{2}} + 1 & \text{ for even } N, \\[.5em]
    \frac{N_q}{2} + 1 &\text{ for odd } N \text{ and } i \leq \frac{N_q}{2}, \\[.5em]     
    \frac{N_q}{2}     &\text{ for odd } N \text{ and } i > \frac{N_q}{2}, \
    \end{cases}
\end{align}
and $\floor{\frac{N_q}{2}}$ denotes integer division.
This assignment of two-point flux computations still accounts for every combination where $j > i$, while leading to 
substantially better load balancing. For even polynomial orders, every thread performs the same number of flux computations, which is optimal. For odd polynomial orders, half the threads do one more iteration, which is still much better than the naive way.    
A modification of this approach is to keep the loop bounds constant at $\frac{N_q}{2} + 1$ and appropriately apply a weight of $1/2$ to antipodal pairs when the polynomial order is odd.
The benefit of this approach is that each thread performs the same number of loop iterations, which can avoid potential issues when synchronizing threads.
The results reported here use the weighting approach, but we include the description of the indexing approach since it illustrates for which pairs the antipodal double counting occurs.
\RevThree{We include Fig.~\ref{fig:thread-mapping} to help explain the thread allocation.}
\input{tikz-figures/thread_mapping}

Figure~\ref{fig:kernel-timings} shows that the combination of 
precomputation of primitive variables and logarithms,
the efficient logarithmic mean formulation, and the exploitation of symmetry accounts for the majority of the total $7.3\times$ speedup. 
\RevOne{These are architecture-independent algorithmic optimizations that can also be applied to CPUs and other accelerators.
Other broadly applicable optimizations shown in Fig.~\ref{fig:kernel-timings} include computing the contravariant flux directly rather than forming and contracting an intermediate $3 \times 5$ flux matrix, directly computing the inverse logarithmic mean, and reusing previously computed flux differences.
The ``micro-optimizations'' step aggregates several small implementation adjustments including reusing previously computed reciprocals and index offsets, moving data loads closer to where the values are used, and rearranging arithmetic expressions.}
\RevOne{Several of the remaining optimizations depend more directly on the GPU execution and memory models.}

The two-point flux involves communication across the entire element.
Therefore, nodes regularly access information from other nodes, making memory efficiency and layout important.
A key aspect of this is placing the state variables, the auxiliary variables, and other regularly accessed quantities such as the precomputed logarithms and differentiation matrix in shared memory.
This gives the whole element quick access to other nodes' values.
Our baseline kernel was launched once per direction, loading shared memory each time.
\RevOne{We refer to combining these three directional kernel launches into a single combined kernel as ``kernel fusion''.
This eliminates redundant shared-memory loads and reduces kernel-launch overhead.}
\RevOne{In addition, we eliminate unnecessary temporary local copies of the state and auxiliary data.
These copies increase per-thread register pressure. Because the available GPU registers are shared among concurrently executing threads, greater per-thread register use can reduce occupancy.}

\RevOne{The distinction between architecture-independent and GPU-specific optimizations is not always absolute.
For example, kernel fusion and elimination of temporary copies can be useful for multiple different architectures.
In the present GPU implementation, however, their primary benefits are reduced kernel-launch overhead, reduced shared-memory traffic, and lower register pressure, all of which are especially beneficial for GPU performance.
}

%% file: tikz-figures/thread_mapping.tex
\begin{figure}[htbp]
    \centering

    \newcommand{\threadgrid}{%
        \foreach \k in {0,...,4} {
            \draw[gray, thin] (\k,0) -- (\k,4);
            \draw[gray, thin] (0,\k) -- (4,\k);
        }
        \foreach \i in {1,...,4} {
            \node at ({\i-0.5},-0.3) {\i};
            \node at (-0.3,{\i-0.5}) {\i};
        }
        \node at (2,-0.9) {Partner node $j$};
        \node[rotate=90] at (-0.9,2) {Thread / node $i$};
        \node[align=center, font=\footnotesize] at (5.1,-0.4)
            {Flux\\evaluations};
    }

    \begin{subfigure}[t]{0.49\textwidth}
        \centering
        \begin{tikzpicture}[x=1cm, y=-1cm, font=\small]
            \foreach \i in {1,...,4}
                \fill[gray!15] ({\i-1},{\i-1}) rectangle (\i,\i);

            \foreach \i/\j/\c in {
                1/2/blue, 1/3/blue, 1/4/blue,
                2/3/orange, 2/4/orange,
                3/4/teal%
            } {
                \fill[\c!25] ({\j-1},{\i-1}) rectangle (\j,\i);
                \node at ({\j-0.5},{\i-0.5}) {$F_{\i\j}$};
            }
            \threadgrid
            \foreach \i/\nflux in {1/3, 2/2, 3/1, 4/0}
                \node at (5.1,{\i-0.5}) {\nflux};
        \end{tikzpicture}
        \caption{Upper-triangular assignment.}
        \label{fig:thread-mapping-triangular}
    \end{subfigure}\hfill
    \begin{subfigure}[t]{0.49\textwidth}
        \centering
        \begin{tikzpicture}[x=1cm, y=-1cm, font=\small]
            \foreach \i in {1,...,4}
                \fill[gray!15] ({\i-1},{\i-1}) rectangle (\i,\i);

            \foreach \i/\j/\c in {
                1/2/blue, 2/3/orange, 3/4/teal, 4/1/violet%
            } {
                \fill[\c!25] ({\j-1},{\i-1}) rectangle (\j,\i);
                \node at ({\j-0.5},{\i-0.5}) {$F_{\i\j}$};
            }
            \foreach \i/\j/\c in {
                1/3/blue, 2/4/orange, 3/1/teal, 4/2/violet%
            } {
                \fill[\c!25] ({\j-1},{\i-1}) rectangle (\j,\i);
                \node at ({\j-0.5},{\i-0.5}) {$\frac{1}{2}F_{\i\j}$};
            }
            \threadgrid
            \foreach \i in {1,...,4}
                \node at (5.1,{\i-0.5}) {2};
        \end{tikzpicture}
        \caption{Weighted cyclic assignment.}
        \label{fig:thread-mapping-cyclic}
    \end{subfigure}

    \caption{\RevThree{Thread assignment for pairwise flux evaluations
    along one coordinate line with $N=3$ and $N_q=4$. Colors identify the
    evaluating thread, and diagonal pairs are omitted.
    (a) The upper-triangular assignment requires three, two, one, and zero
    evaluations from threads 1--4, respectively.
    (b) The weighted cyclic assignment used in this work performs two
    evaluations per thread.
    The antipodal pairs $\{1,3\}$ and $\{2,4\}$ are evaluated in both
    orientations, with weight $1/2$ on each evaluation so that each unordered
    pair has total weight one.}}
    \label{fig:thread-mapping}
\end{figure}

%% file: sections/performance_analysis.tex
\subsection{Single-GPU}
\label{sec:single-gpu}

\begin{figure}[htbp]
    \centering
    \begin{subfigure}[t]{0.49\textwidth}
        \centering
        \includegraphics[width=\textwidth]{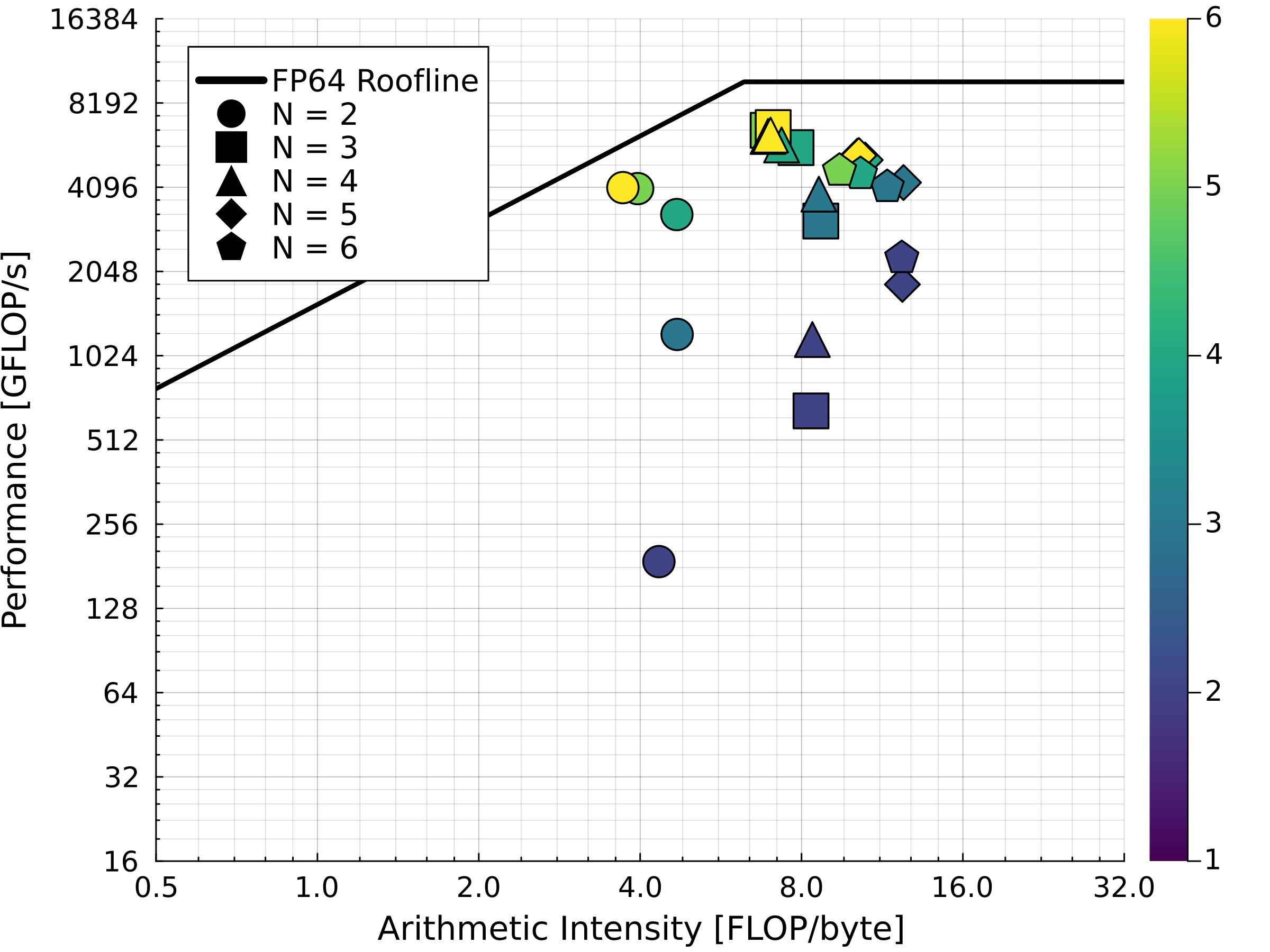}
        \caption{All refinement levels.}
        \label{fig:volume-knl-roofline-all}
    \end{subfigure}\hfill
    \begin{subfigure}[t]{0.49\textwidth}
        \centering
        \includegraphics[width=\textwidth]{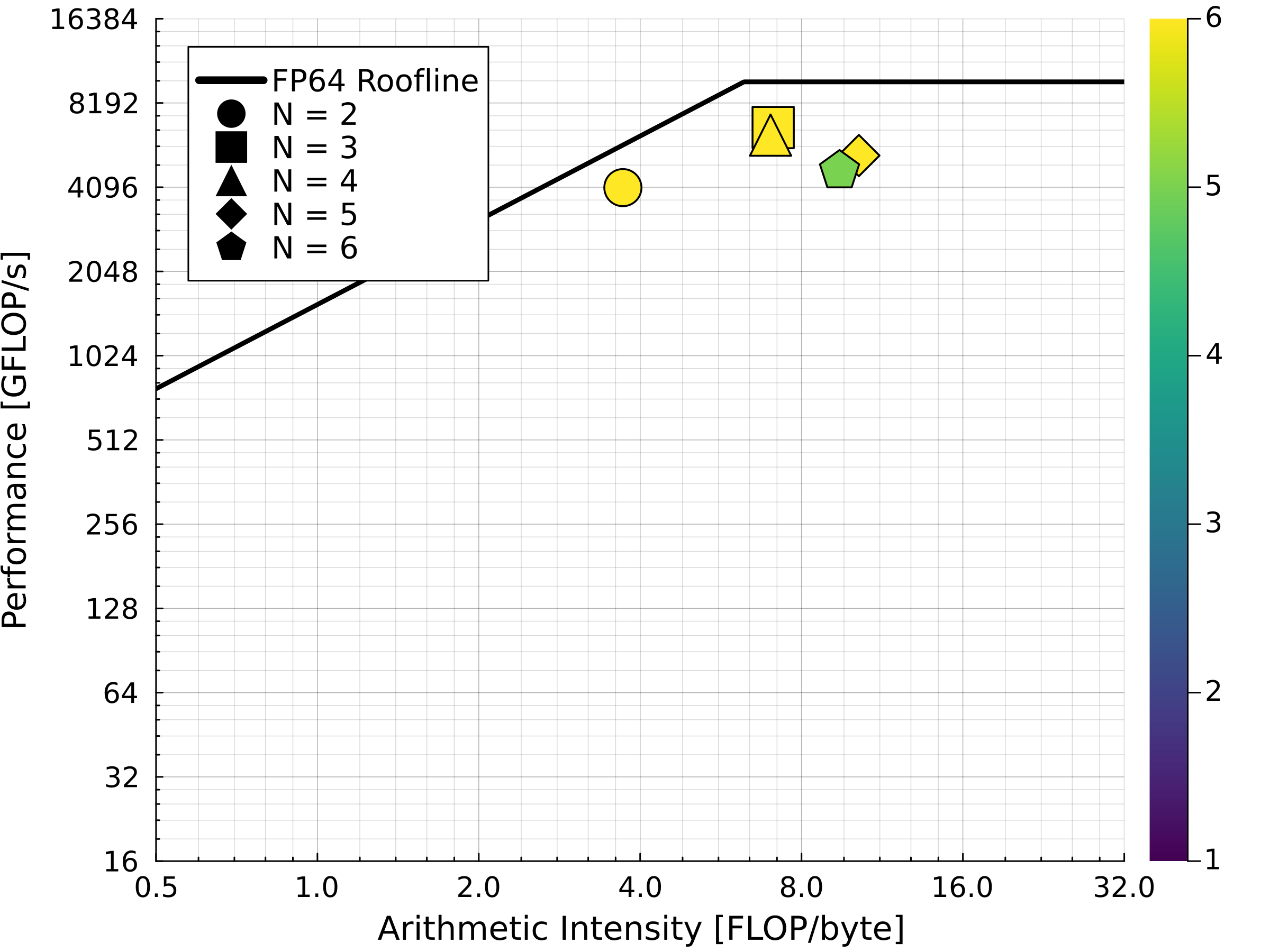}
        \caption{Most-saturated run for each $N$ ($L=6$ for $N \le 5$; $L=5$ for $N=6$, which was out of memory at $L=6$).}
        \label{fig:volume-knl-roofline-sat}
    \end{subfigure}
    \caption{Roofline plots of the optimized volume kernel on an NVIDIA A100 SXM4 40GB (\RevOne{peak 64-bit performance}: 9{,}746 GFLOP/s, peak bandwidth: 1{,}560 GB/s) \RevOne{using the rising thermal bubble test from Sec.~\ref{sec:rising-bubble-case}}. Marker shape denotes polynomial order $N$ and color denotes refinement level $L$ (dark to light: coarse to fine). \textbf{(\subref{fig:volume-knl-roofline-all})} shows every $(N,L)$ configuration; \textbf{(\subref{fig:volume-knl-roofline-sat})} restricts to the most refined mesh available for each polynomial order, isolating the saturated-GPU regime. 
    }
    \label{fig:volume-knl-roofline}
\end{figure}

Our optimized volume kernel is in the compute-bound regime for polynomial order $N \geq 3$, achieving nearly 70\% of the roofline ceiling.
We can show this by analyzing a roofline plot of the kernel, a performance model introduced in \cite{williams_roofline_2009}.
We plot the result in Fig.~\ref{fig:volume-knl-roofline} with the results for all refinement levels on the left and the most saturated results on the right.
The horizontal axis is arithmetic intensity, which is a ratio of floating-point operations to bytes transferred.
The vertical axis is performance, which is the number of floating-point operations per second.
As we increase the polynomial order, denoted by marker shape, our data points move to the right since there is more work to be done per piece of data on the hardware.
As we increase the refinement level, denoted by marker color, we increase the amount of parallel work, pushing performance upward.
We observe a slight leftward shift at higher refinement levels due to less optimal cache usage for larger problem sizes.
We also see 
decreasing achieved throughput as polynomial order increases
beyond $N = 3$, which is attributable to larger thread block sizes reducing occupancy.

DG methods are typically run at high orders ($N \geq 3$)  and the plot shows that all such cases are compute-bound.
From this plot, we can infer that (i) the volume kernel is well-suited to GPU parallelization and makes full use of the computational resources since it is not memory bandwidth-limited 
at typical polynomial orders and (ii) our code can take advantage of rapidly evolving GPU hardware since compute-bound kernels benefit from FLOP/s increases. 

\begin{figure}[htbp]
    \centering
    \begin{subfigure}[t]{0.48\textwidth}
        \centering
        \includegraphics[width=\textwidth]{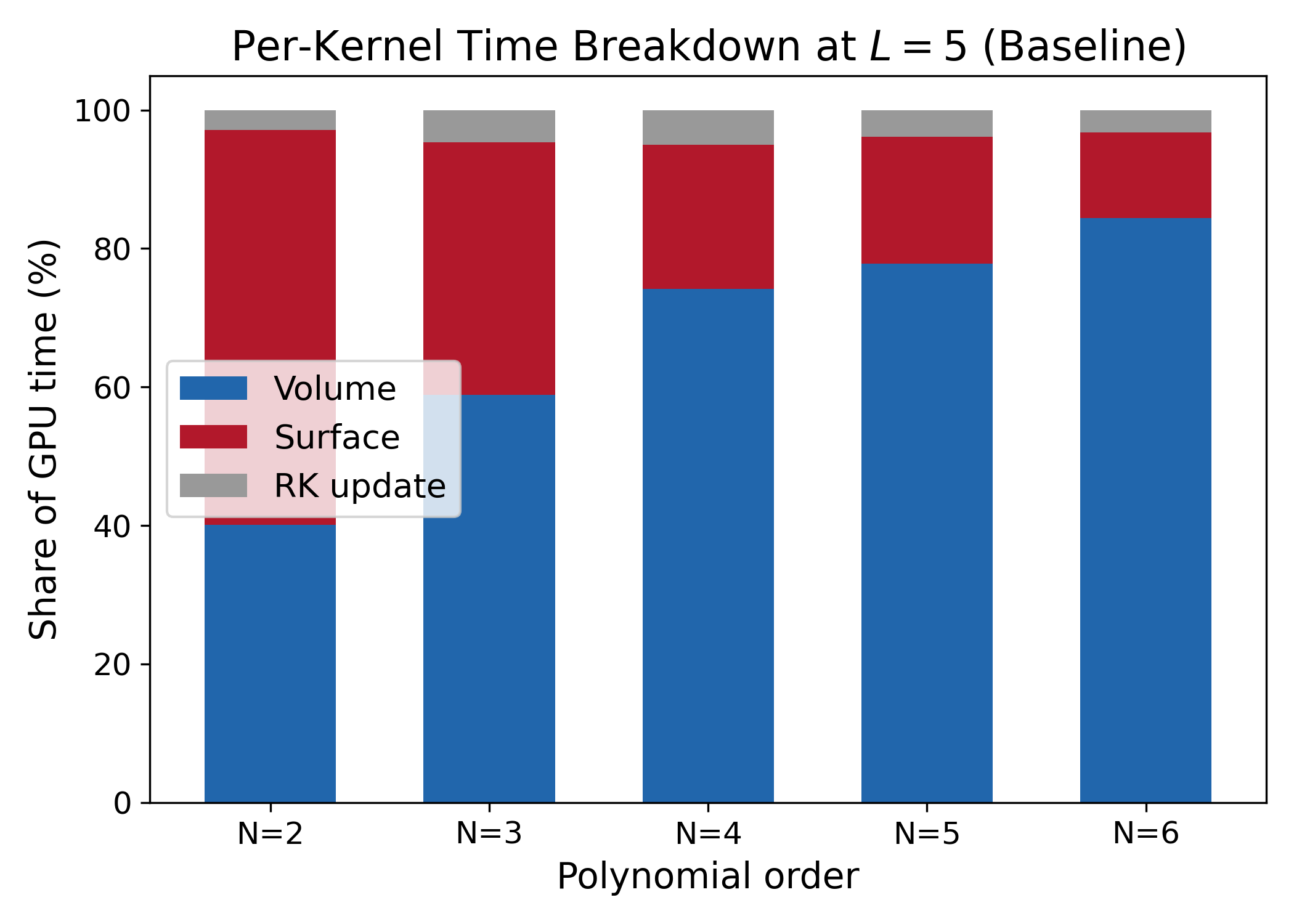}
        \caption{Per-kernel share of GPU time at $L=5$ with baseline volume kernel.}
        \label{fig:gpu-time-percents-breakdown-baseline}
    \end{subfigure}\hfill
    \begin{subfigure}[t]{0.48\textwidth}
        \centering
        \includegraphics[width=\textwidth]{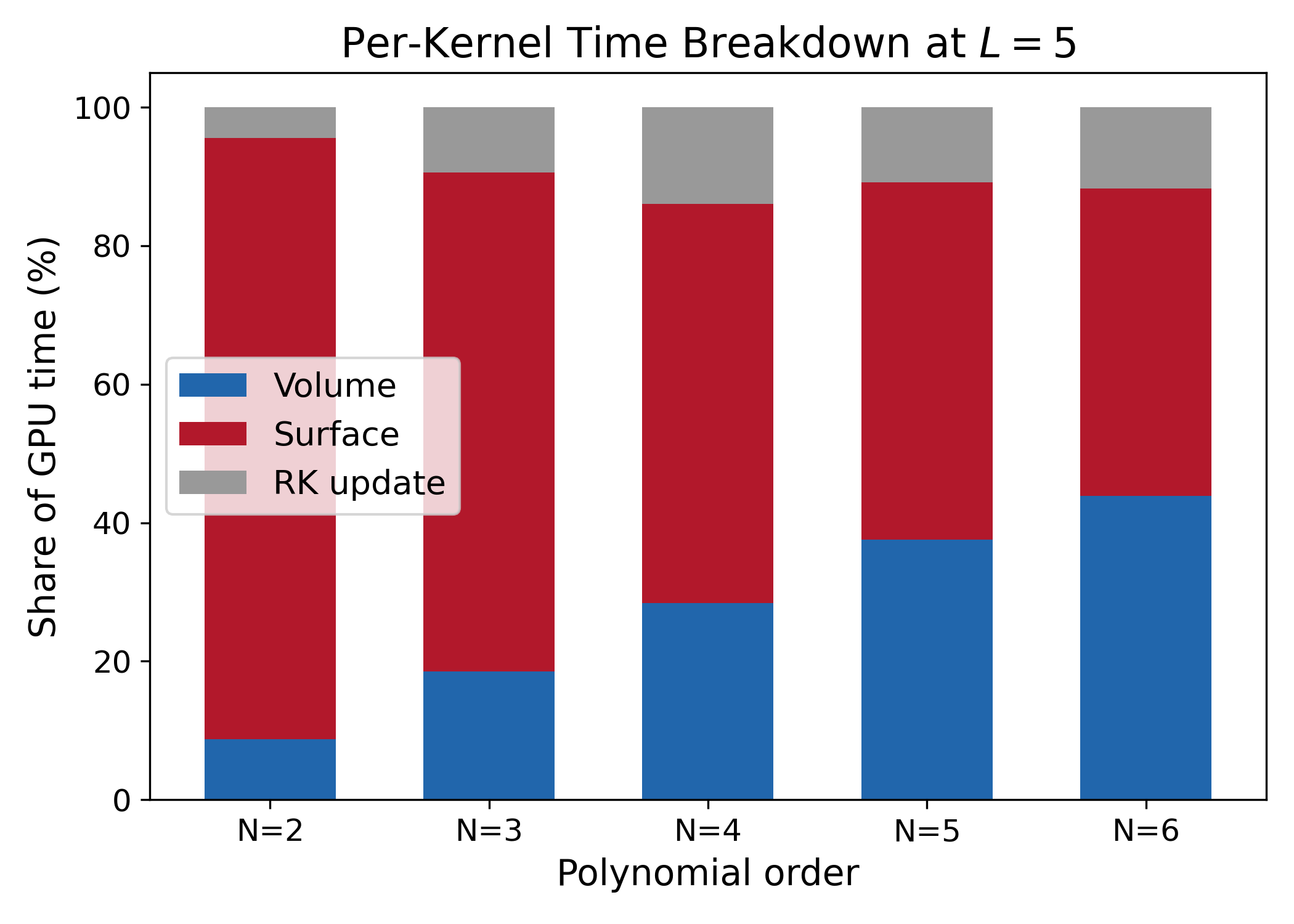}
        \caption{Per-kernel share of GPU time at $L=5$ with optimized volume kernel.
        }
        \label{fig:gpu-time-percents-breakdown-optimized}
    \end{subfigure}
    \caption{Percent of total GPU time for the three main kernels at refinement level $L=5$ \RevOne{using the rising-bubble test case from
    Sec.~\ref{sec:rising-bubble-case} as the underlying simulation}.
    \textbf{(\subref{fig:gpu-time-percents-breakdown-baseline})} The unoptimized volume kernel dominates for \RevOne{polynomial order} $N\geq 3$.
    \textbf{(\subref{fig:gpu-time-percents-breakdown-optimized})} After optimization, the surface kernel becomes the dominant cost across the tested polynomial orders.
    }
    \label{fig:gpu-time-percents}
\end{figure}

A complete stage of one time step involves three parts: the volume kernel, the surface kernel, and the Runge-Kutta (RK) time-stepping kernel; we have only discussed the volume kernel so far.
Figure~\ref{fig:gpu-time-percents} shows the per-kernel share of GPU time before and after our volume-kernel optimizations.
In panel~(\subref{fig:gpu-time-percents-breakdown-baseline}), the volume kernel is the bottleneck for $N \geq 3$ and dominates total runtime; in panel~(\subref{fig:gpu-time-percents-breakdown-optimized}), after optimization, the surface kernel becomes the dominant cost across all tested polynomial orders. The optimized volume kernel therefore addresses the cost issue of ESDG, and future optimizations will likely need to target the surface kernel.
\RevOne{Surface-flux evaluation has low arithmetic intensity.
At each face degree of freedom, state and auxiliary variables from both sides of the interface are loaded along with face-mapping and local geometric data, but are used for only one flux evaluation.
Performance is therefore limited primarily by memory bandwidth.
The repeated pairwise evaluations and substantial data reuse targeted by the volume-kernel optimizations are absent, limiting the applicability of those techniques.
Further improvements must instead reduce data transfer or improve memory-access patterns.}
\begin{figure}[htbp]
    \centering
    \includegraphics[width=0.5\textwidth]{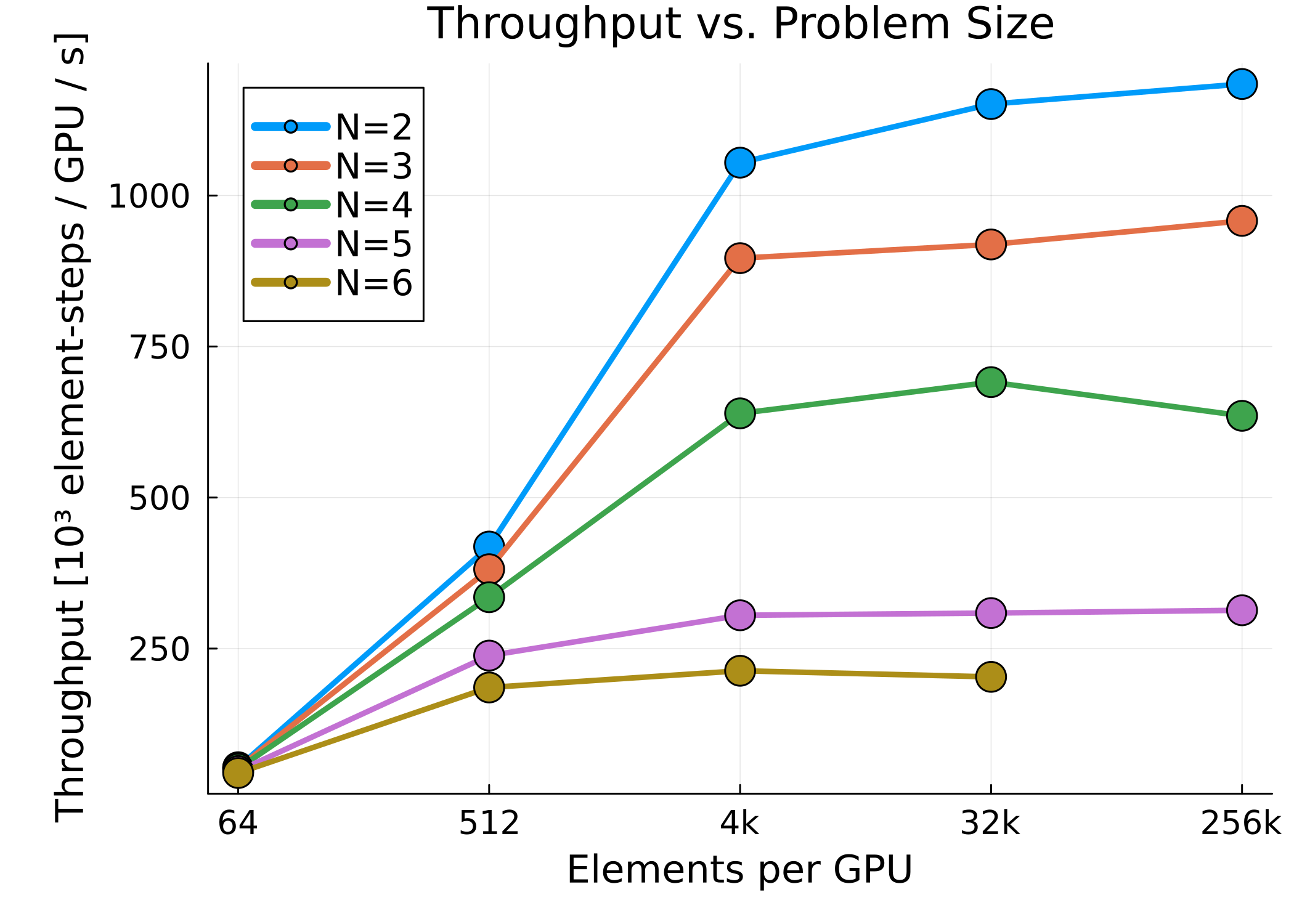}
    \caption{Single-GPU throughput in thousands of element--time-steps per second as the number of elements grows (\RevOne{a base mesh consisting of a single coarse element ($K=1$) refined $L$ times, giving an element count $= 2^{3L}$}).  
    Each curve corresponds to one polynomial order; the flattening of each curve marks the element count at which the GPU is saturated for that order. \RevOne{Although we used the rising thermal bubble described in Sec.~\ref{sec:rising-bubble-case} as the underlying simulation, the results should be test-case agnostic}.
    }
    \label{fig:throughput}
\end{figure}

Inspired by \cite{Bertagna2019}, we show the throughput as a function of the number of elements in Fig.~\ref{fig:throughput}; on the horizontal axis we vary the number of three-dimensional elements on a single GPU while on the vertical axis we report on the number of elements evaluated in one time step per wall-clock second (in thousands). The higher values for the lowest-degree polynomial ($N=2$) show that even at 256k elements we have not yet fully saturated the device. In contrast, for the highest-degree polynomial ($N=6$) we are able to fully saturate the device at 4k elements. 
This difference in the number of elements needed for saturation can be attributed to elements at higher polynomial orders containing more degrees of freedom and, therefore, having more work per degree of freedom (i.e., higher arithmetic intensity). 
For example, for a problem with 128k elements and polynomial order $N=4$, the optimal number of GPUs to deploy would be 4 (resulting in 32k elements per GPU) since each GPU is fully saturated at this point.
If we used more GPUs, the results would shift leftward along the curve, thereby underutilizing each device. 

\RevThree{As a final consideration for the single-GPU case, we compare the performance of our ESDG kernels with our best previously published results for standard DG kernels \cite{Abdi2017a,Abdi2017b}. 
Within the roofline model, the crossover arithmetic intensity separating the memory-bandwidth-bound and compute-bound regimes is
\[
I_\text{crossover} = \frac{ P_\text{peak} } { {B_\text{peak}} }
\]
where $P_\text{peak}$ is the peak floating-point performance in GFLOP/s, $B_\text{peak}$ is the peak memory bandwidth in GB/s, 
and \(I_{\text{crossover}}\) is expressed in GFLOP/GB.
A kernel with arithmetic intensity $I < I_\text{crossover}$ is expected to be memory-bandwidth-bound, whereas one with $I > I_\text{crossover}$ is expected to be compute-bound.
For an NVIDIA A100 GPU, the crossover arithmetic intensity is
\[
I^\text{A100}_\text{crossover} = \frac{ 9746 \, \text{GFLOP/s} } { {1560 \text{GB/s} } } = 6.25 \, \text{GFLOP/GB}.
\]
For the NVIDIA GTX Titan Black GPU used in \cite{Abdi2017a,Abdi2017b}, it is
\[
I^\text{GTX}_\text{crossover} = \frac{ 1707 \, \text{GFLOP/s} } { {334 \text{GB/s} } } = 5.11 \, \text{GFLOP/GB}.
\]
Figure~\ref{fig:volume-knl-roofline} shows that nearly all ESDG volume kernels have arithmetic intensity greater than 4 GFLOP/GB, with the highest intensity kernels exceeding 7  GFLOP/GB. In contrast, the standard DG kernels shown in Fig.\ 5(d) of \cite{Abdi2017a} have arithmetic intensity of approximately 3 GFLOP/GB, corresponding to 
a 2.5$\times$ increase for the ESDG kernels.  Although the ESDG kernels perform substantially more computation than the standard DG kernels, they also reuse data more effectively, thereby increasing their arithmetic intensity and making better use of the available computational resources.  Consequently, the standard DG volume kernel achieves approximately 17\% of peak floating-point performance, whereas the ESDG volume kernel achieves 67\%.}


\subsection{Multi-GPU}
\label{sec:multi-gpu}

The preceding section focused on performance analysis on one GPU; we now move to  multi-GPU strong and weak scaling results. Before discussing the scaling results, we first describe the Message Passing Interface (MPI) strategy. 
We assign each MPI rank a single GPU and distribute elements to these GPUs using \texttt{p4est}; more information on the 
distribution method based on a Morton space-filling curve
can be found in \cite{BursteddeWilcoxGhattas11}.
We implement the overlapping communication strategy of \cite{Kelly2012}.
This entails posting non-blocking sends and receives for the boundary data at the beginning of each Runge-Kutta stage so that the communication overlaps with both the volume flux and the intra-processor surface integrals.
The algorithm then waits for the communication to finish and computes the inter-processor surface integrals.
\RevOne{The multi-GPU simulations are run with one MPI rank per GPU and CUDA-aware MPI with GPUDirect support.
This enables direct GPU-to-GPU transfers within a node and communication between the network interface and GPUs across nodes.
Enabling this sort of communication requires few code-level changes but does require appropriate support and configuration on the HPC system.
\RevTwo{For all scaling experiments, runs with up to four GPUs use a single node, while larger simulations use four GPUs per node.}
GPU timings use CUDA-event durations for each timestep.
For each timestep, we take the maximum CUDA-event duration across MPI ranks.
We sum these maxima for total times and take the median for weak-scaling plots.
All recorded timesteps are included.}

\begin{figure}[htbp]
    \centering
    \begin{subfigure}[t]{0.48\textwidth}
        \centering
        \includegraphics[width=\textwidth]{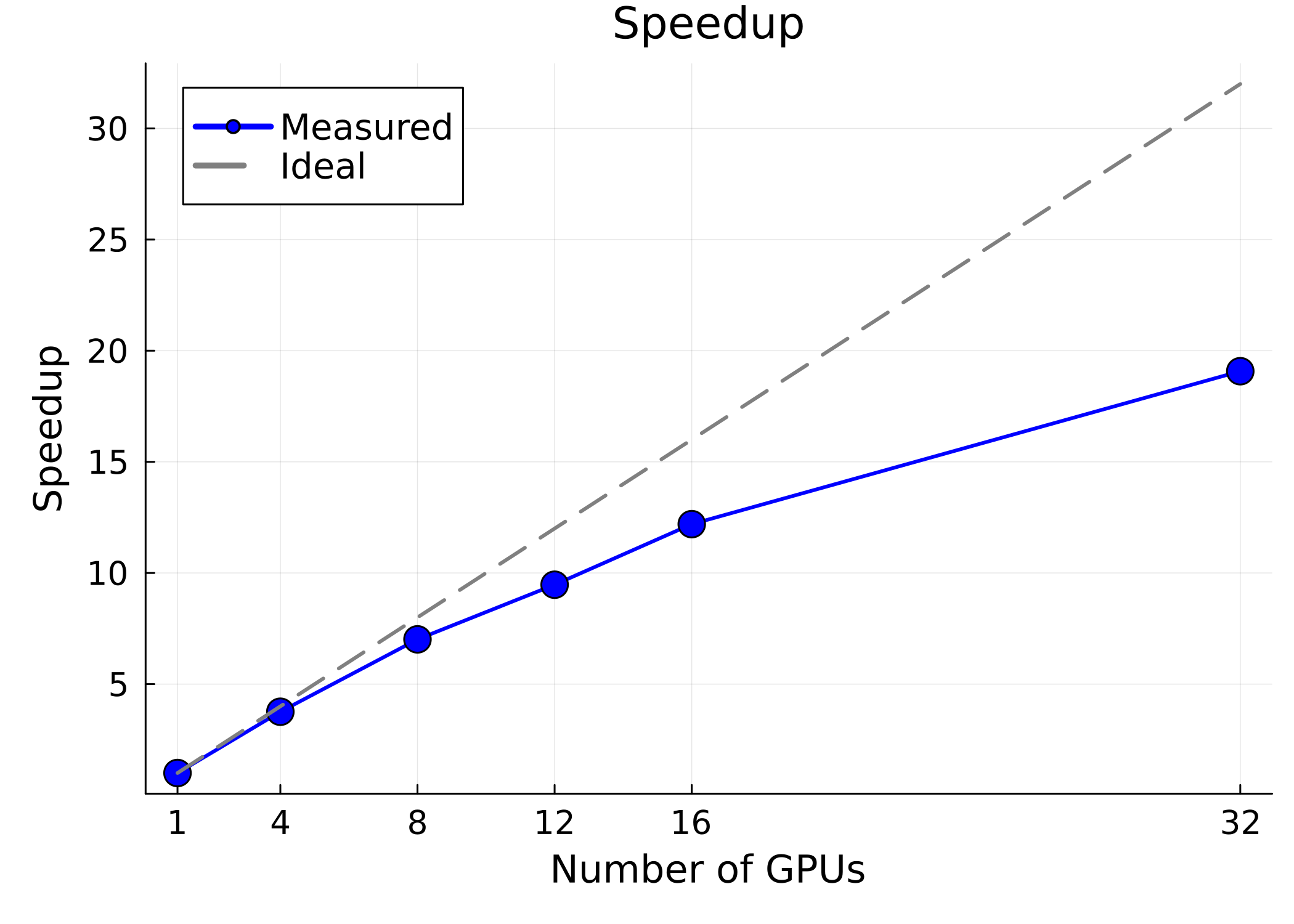}
        \caption{Strong scaling.}
        \label{fig:scaling-results-strong}
    \end{subfigure}\hfill
    \begin{subfigure}[t]{0.48\textwidth}
        \centering
        \includegraphics[width=\textwidth]{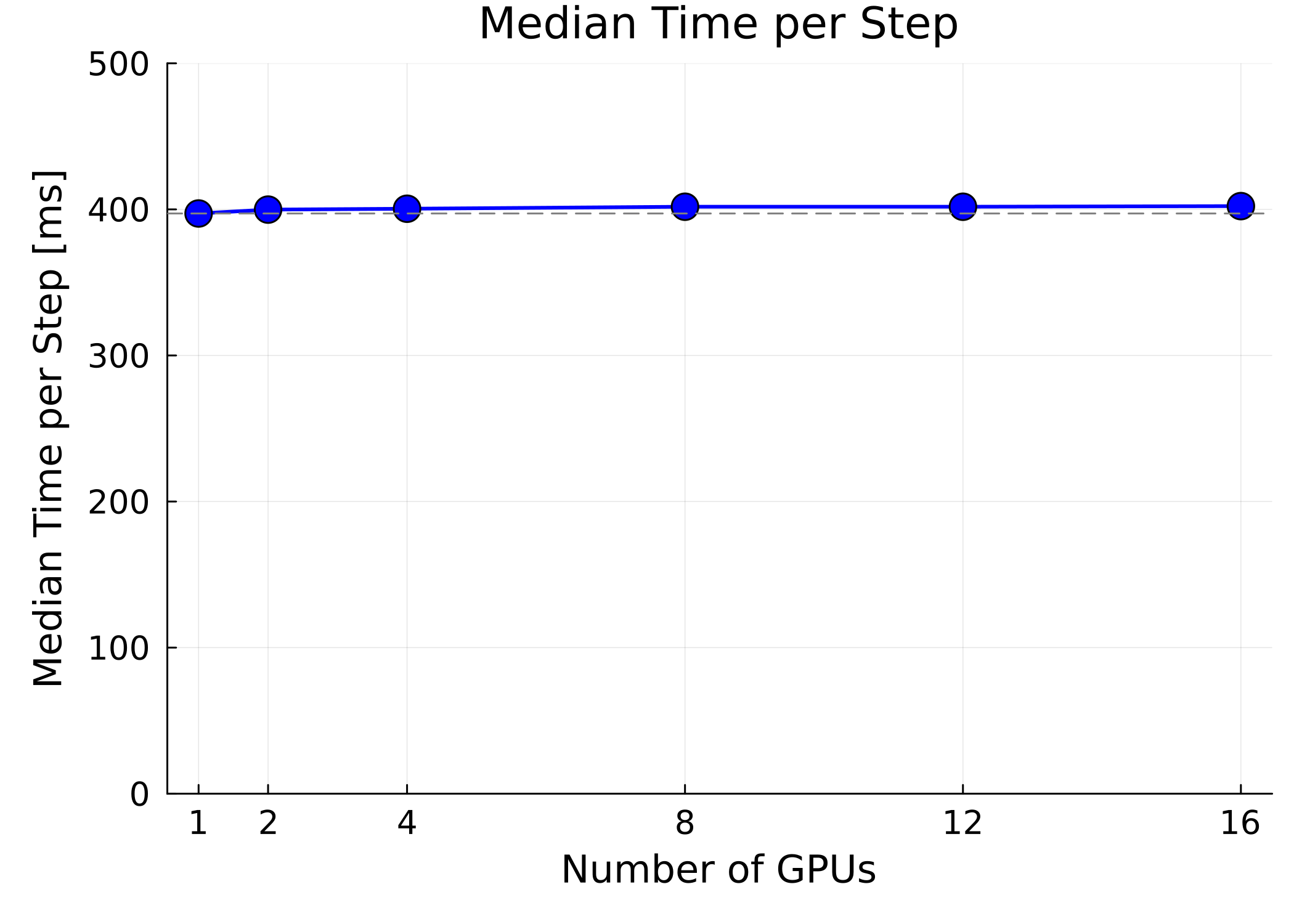}
        \caption{Weak scaling.}
        \label{fig:scaling-results-weak}
    \end{subfigure}
    \caption{Multi-GPU scaling of the rising-bubble benchmark at \RevOne{polynomial order} $N{=}4$ and \RevOne{refinement level} $L{=}6$ on NVIDIA A100 GPUs \RevTwo{with 4 GPUs per node}. \textbf{(\subref{fig:scaling-results-strong})} strong-scaling speedup on 1--32 GPUs for a fixed 32.8M-DOF problem, against ideal linear speedup. \textbf{(\subref{fig:scaling-results-weak})} weak-scaling median time per step as the element count grows with the GPU count, \RevThree{with each device maintaining the same number of elements (optimal graph partitioning)}.}
    \label{fig:scaling-results}
\end{figure}

Strong scaling fixes the problem size and increases the processor count.
For this test, we run the rising-bubble test case with fourth-order polynomials and 6 refinement levels ($2^6 = 64$ elements in each direction) for 1000 time steps.
We show the results in the left panel of Fig.~\ref{fig:scaling-results}.
The left panel shows that speedup falls below ideal as GPU count increases.
\RevOne{As the GPU count increases, the workload per GPU decreases, leaving less computational work to overlap with the boundary-data exchange.
Although the boundary-data messages also become smaller, the communication latency inherent in sending messages of any size does not decrease proportionally.
Therefore, an increasing fraction of this communication latency is exposed, contributing to the non-ideal strong scaling.}
\RevOne{The corresponding \texttt{p4est} data shows that each rank receives the same number of elements, indicating negligible computational load imbalance.
However, at 16 and 32 GPUs, half of the ranks have send patterns containing approximately $34\%$ and $26\%$ more interface nodal indices, respectively, than the other half.
These results indicate that the departure from ideal strong scaling reflects reduced per-GPU utilization and increasingly exposed communication costs, including communication-load imbalance at larger GPU counts, rather than imbalance in the computational workload.}
At 32 GPUs, parallel efficiency drops to 60\%, suggesting that a better use of the hardware is to match problem size to processor count, which we explore below.

For weak scaling, we increase both problem size and number of processors simultaneously.
For this test we once again used polynomial order $N=4$ and refinement level $L=6$, but this time we increased the number of coarse elements to coincide with the number of processors---the results are shown in the right panel of Fig.~\ref{fig:scaling-results}.
\RevThree{For example, the run with P processors has P coarse elements of polynomial order $N=4$ and refinement level $L=6$.}
When we increase to 8 processors, we now have the same setup but with 8 elements.
This is representative of a real-world use case in which the number of processors is selected based on the size of the problem.
We see from our results that our code maintains approximately the same time per step as the GPU count and problem size increase.
This is an ideal result that shows the time to complete a fixed number of numerical steps for a large problem can be approximately the same as that of a smaller problem as long as we properly allocate computational resources.

\begin{figure}[htbp]
    \centering
    \begin{subfigure}[t]{0.48\textwidth}
        \centering
        \includegraphics[width=\textwidth]{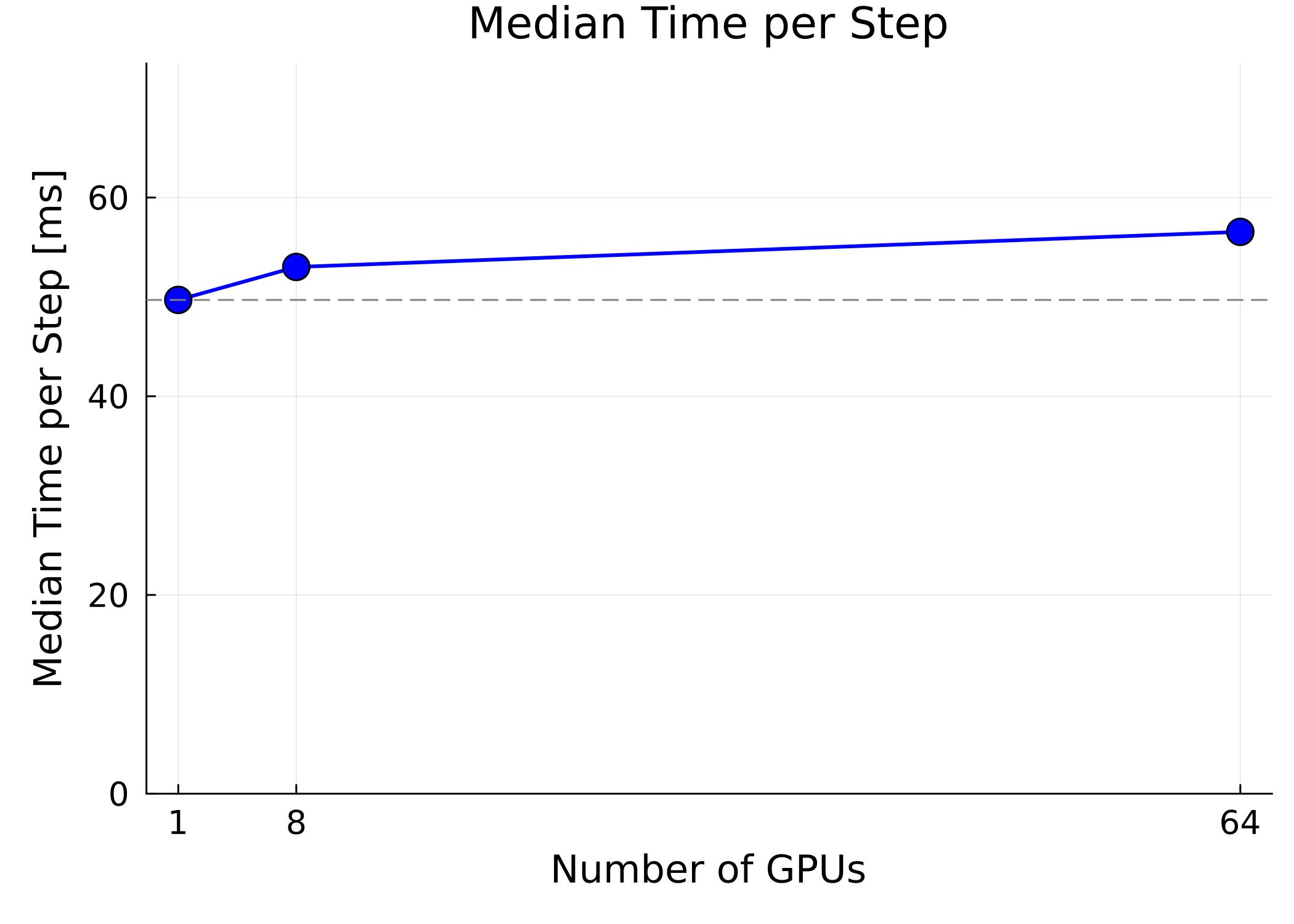}
        \caption{Includes $P{=}1$ baseline.}
        \label{fig:weak-scaling-sfc-with-single}
    \end{subfigure}\hfill
    \begin{subfigure}[t]{0.48\textwidth}
        \centering
        \includegraphics[width=\textwidth]{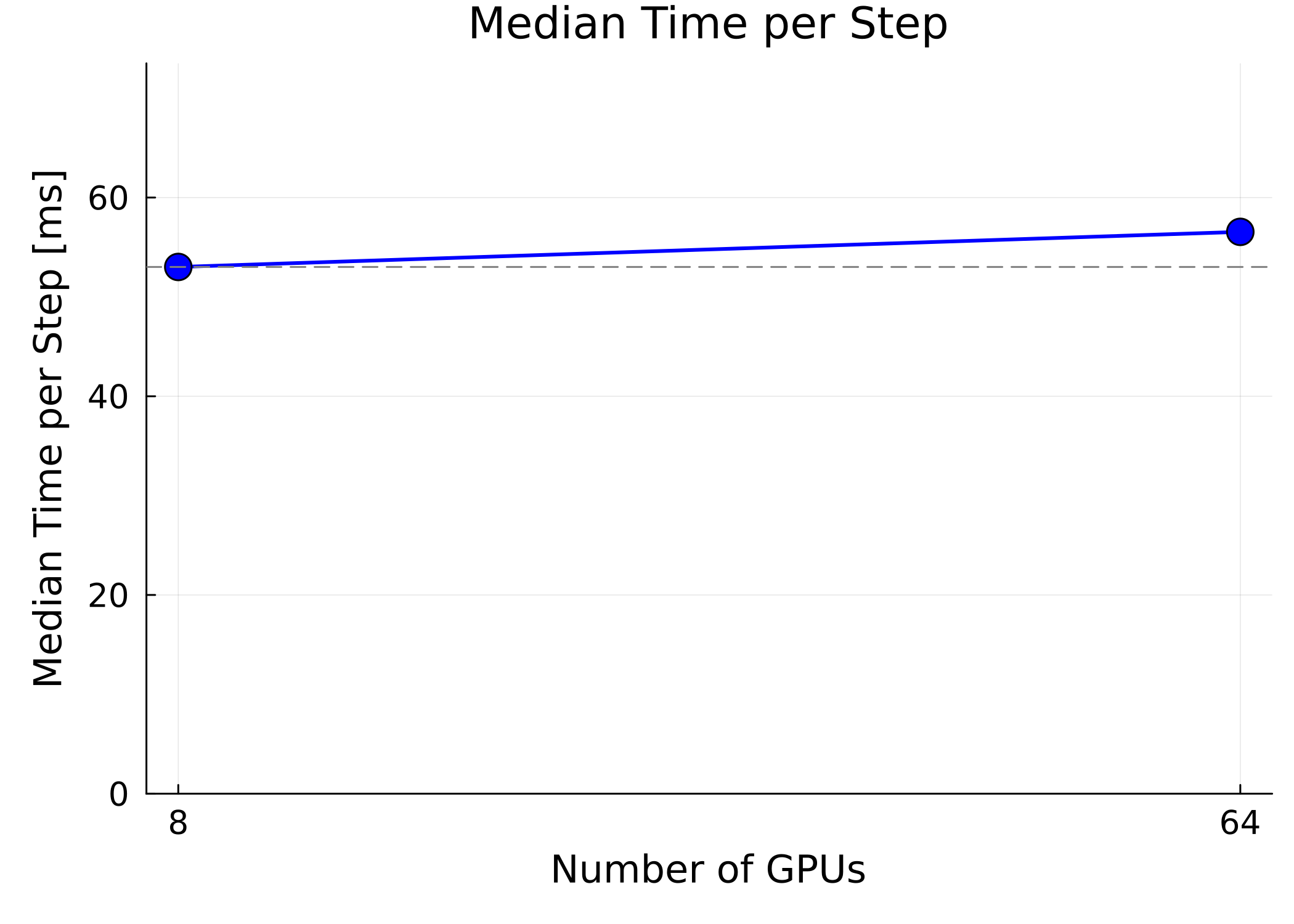}
        \caption{Multi-GPU only.}
        \label{fig:weak-scaling-sfc-no-single}
    \end{subfigure}
    \caption{Weak-scaling median time per step for the rising-bubble benchmark where the GPU count grows with the refinement level (a single coarse element is subdivided so each GPU holds identical work: \RevOne{refinement level} $L{=}6$ at 8 GPUs, $L{=}7$ at 64 GPUs). \textbf{(\subref{fig:weak-scaling-sfc-with-single})} includes the \RevOne{processor count} $P{=}1$ baseline, which has no inter-GPU communication. \textbf{(\subref{fig:weak-scaling-sfc-no-single})} baseline shifted to $P{=}8$ so only multi-GPU effects are visible.}
    \label{fig:weak-scaling-sfc}
\end{figure}

\RevThree{The preceding weak scaling test is "ideal" in that the grid refinements are each in separate coarse elements and the number of coarse elements equals the number of processors. 
This guarantees a minimum number of edge cuts per device, which translates into fewer neighbors.
To stress our code with a less ideal scenario, we use one coarse element and increase the refinement level and processor count to keep the problem size per processor constant---in contrast to the ideal case, each device will have more neighbors to communicate with.}
These results are shown in the left panel of Fig.~\ref{fig:weak-scaling-sfc}.

The drawback of this test is that the exponential increase in the number of required GPUs to scale with the refinement limits the data we can gather since we do not have access to 512 GPUs.
For this reason, we include the single-GPU case in the left panel of Fig.~\ref{fig:weak-scaling-sfc}, but clearly the single GPU will not have communication overhead and we see a large drop in efficiency from 1 to 8 GPUs.
We remove this single-GPU data point in the right panel of Fig.~\ref{fig:weak-scaling-sfc} and see that the efficiency drop between multi-GPU runs is not as significant.

\subsection{GPU versus CPU}
\RevThree{We next compare the performance of the GPU and CPU codes for the rising thermal bubble using polynomial order \(N=4\) and refinement level \(L=6\) over 1,000 time steps.} For the CPU code, we use the \RevTwo{Nonhydrostatic Unified Model of the Atmosphere (NUMA)}  \cite{Giraldo2013,Giraldo2024}. For the comparison, we use exactly the same grid size, time step, and numerical methods (ESDG with the explicit \RevTwo{fourth-order fourteen-stage Low-Storage Runge-Kutta (LSRK)} time-stepping). Although the GPU and CPU codes are rather different (the GPU code is written in Julia~\cite{bezanson2017julia} using CUDA.jl~\cite{besard2018juliagpu} with MPI.jl~\cite{byrne2021mpi}  and the CPU code is written in modern Fortran using MPI), both codes use \texttt{p4est}~\cite{BursteddeWilcoxGhattas11} to generate meshes in parallel and have been designed to run optimally on their specific target hardware platform (e.g., see NUMA results on CPUs in \cite{Mueller2018}).\footnote{\RevTwo{In that study, $N=3$ was shown to be optimal on the CPU but here we use $N=4$ for both the GPU and CPU because it represents the highest order we can use without running into indexing issues as described in Fig.~\ref{fig:thread-mapping}, while still permitting  $L=6$ grid refinement on a single A100 GPU.}} For the comparison, we run the codes on the Delta computer at the National Center for Supercomputing Applications, which has 400 NVIDIA A100 GPUs and 256 AMD EPYC 7763 Milan nodes; in the results below we report the CPU result for 2 sockets per node, so effectively Delta has 128 of the AMD Milan dual-socket nodes.  The thermal design power (TDP) of the A100 is 400 watts, while that for an AMD (7763) socket (with 64 cores) is 280 watts (we multiply this TDP by 2 since each node has two sockets).\footnote{\RevThree{
We estimate energy consumption from TDP, assuming that both the GPU and CPU operate at full utilization. However, measurements from NVIDIA-SMI and AMDuProf indicate that, on average, the GPU draws 77\% of its TDP, whereas the CPU draws 95\% of its rated TDP.}}
We summarize the comparison between GPU and CPU hardware in Table~\ref{table:gpu_vs_cpu}, where energy consumption is calculated as follows
\RevTwo{
\be
\text{Energy Consumption (MJ)} =\text{No.\ Devices} \times \text{Time-to-solution} \times \text{TDP}/10^6
\label{eq:energy_consumption}
\ee
where $No.\ Devices$ denotes the number of devices (CPU or GPU, defined as the device  with the TDP values defined above), the time-to-solution is given in seconds, and the TDP is given in watts (joules per second).
}
$R_{\mathrm{time}}$ denotes the ratio of the CPU time to the GPU time for the same number of nodes (as described previously). $R_{\mathrm{energy}}$ denotes the ratio of CPU energy consumption to GPU energy consumption. 
\begin{table}
\begin{center}
\begin{tabular}{ c  r  r  r  r  r }
        \toprule
 Hardware & No. Devices & Time-to-Solution (s) & Energy Consumption (MJ) & $R_{\mathrm{time}}$ & $R_{\mathrm{energy}}$ \\ 
        \midrule
 A100  & 4 & 106.88 & 0.17 & 11.61 & 16.26 \\
 Milan & 4 & 1241 & 2.78 & --- & --- \\ \midrule
 A100  & 8 & 57.25 & 0.18 & 10.36 & 14.50 \\
 Milan & 8 & 593 & 2.66 & --- & --- \\ \midrule
 A100  & 16 & 32.92 & 0.21 & 9.78 & 13.69 \\
 Milan & 16 & 322 & 2.89 & --- & --- \\
 \bottomrule
\end{tabular}
\end{center}
\caption{Comparison of time-to-solution, energy consumption, and the ratios between the A100 (GPU) and the AMD Milan (CPU) for the rising thermal bubble.}
\label{table:gpu_vs_cpu}
\end{table}
$R_{\mathrm{time}}$ shows that the GPU is approximately $10\times$ faster than the CPU. To be more specific, this result shows that to achieve the same time-to-solution as one NVIDIA A100 GPU, we would need 10 dual-socket AMD Milan nodes (1280 AMD cores). However, since energy consumption is an important consideration, we also compute $R_{\mathrm{energy}}$, which shows that the GPU is at least $13\times$ more efficient than the CPU.

%% file: sections/simulation_results.tex
\section{Simulation Results}
\label{sec:simulations}

In this section, we demonstrate that the GPU-efficient algorithms described in Sec.~\ref{sec:performance} yield good solutions for geophysical fluid dynamics problems of interest. Below, we show results for a three-dimensional discontinuous rising thermal bubble and for a three-dimensional baroclinic instability in a channel.  For both test cases, we show results in 64-bit and 32-bit floating-point precision.  For the baroclinic instability test, we also compare the ESDG solution against our most mature nonhydrostatic atmospheric model.

\subsection{Rising Bubble}
\label{sec:rising-bubble-case}
\begin{figure}[htbp]
    \centering
    \includegraphics[width=0.8\textwidth]{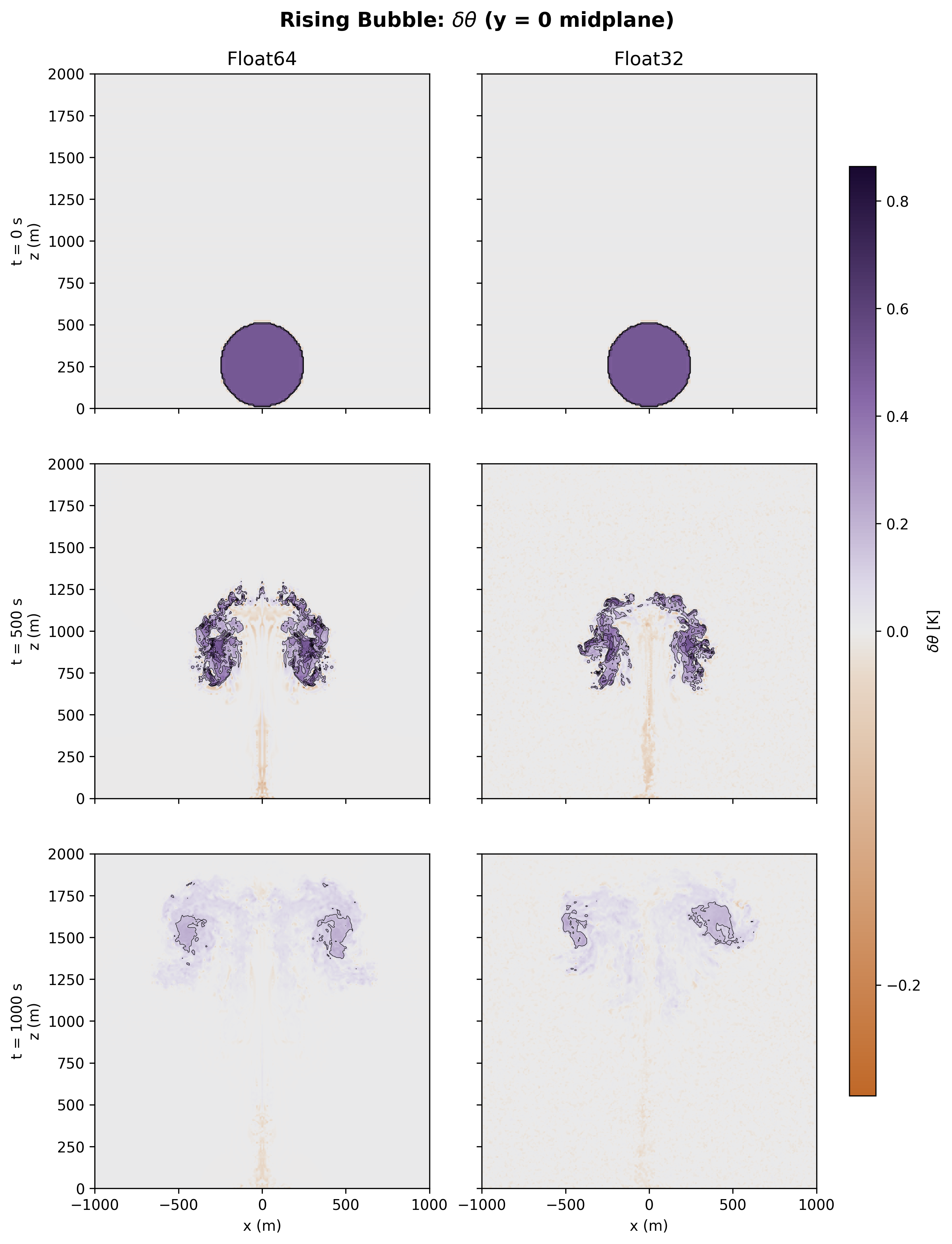}
    \caption{Slices of simulation results for potential temperature perturbation at $y=0$ with polynomial order $N=4$, refinement level $L=6$, and one coarse element in 64-bit precision (left) and 32-bit precision (right). This simulation uses the sharp initial perturbation defined in (\ref{eq:sharp-rb-ic}) for the rising thermal bubble.}
    \label{fig:rising-bubble-floatcomp-sharp}
\end{figure}

We run a full 3D simulation of a rising thermal bubble \cite{Waruszewski2022}.
Our domain is $x \in [-1000, 1000] ~\text{m}$ (periodic), $y \in [-1000,1000]~ \text{m}$ (periodic), and $z \in [0, 2000]~ \text{m}$ (reflecting).
We start the bubble in two different ways.
The first, shown in Fig.~\ref{fig:rising-bubble-floatcomp-sharp}, is a sharp discontinuity defined by: 
\begin{align}
   \delta \theta = 
   \begin{cases}
        \frac{1}{2} & r \leq r_c\\
        0 & \text{otherwise}
   \end{cases} \label{eq:sharp-rb-ic}
\end{align}
where $\delta \theta$ is the potential temperature perturbation, $r = \abs{\bm{x} - \bm{x}_c}$, $r_c = 250~ \text{m}$, and $\bm{x}_c = (0,0,260)~ \text{m}$.
This places a spherical bubble of radius $250~\text{m}$ centered at $(0,0,260)~\text{m}$.
This is a sharp profile that 
jumps abruptly from zero to the full perturbation amplitude.
The second configuration \cite{Giraldo2008} keeps everything the same, but smooths the initial perturbation as follows:
\begin{align}
   \delta \theta = 
   \begin{cases}
        \frac{1}{2} \frac{(1 + \cos(\pi\frac{r}{r_c}))}{2} & r \leq r_c\\
        0 & \text{otherwise}.
   \end{cases}\label{eq:smooth-rb-ic} 
\end{align}
We run both configurations in 64- and 32-bit precision with polynomial order $N =4$ and refinement level $L=6$ on one coarse element, resulting in a total DOF count of 32,768,000 ($N_e=(2^6)^3=262,144$ and $N_p=N_e(N+1)^3=32,768,000$). 
We run this simulation on 8 GPUs with a Courant number ($C$) of 0.5 and a resulting time step of $\Delta t \approx 8.18\times10^{-3}$~s.
We present a slice of the potential temperature results at $y = 0$ for three different time steps.
We see that the 64-bit precision results mostly maintain large-scale symmetry while the 32-bit precision results diverge slightly over time.
Due to the minimal amount of numerical dissipation, the solution may lose symmetry from non-symmetric rounding errors. (For example, the symmetry in Fig.~3 of \cite{Waruszewski2022} was obtained by ensuring rounding errors were made symmetrically.)
Despite these differences in symmetry, the practical benefit of 32-bit precision is substantial: the timing results, summarized in Table~\ref{tab:rb-timing}, show that 32-bit precision consistently yields a $2\times$ speedup across both initial conditions. 
Changing the precision type is as easy as passing an argument to the run script and requires no other configuration.
\begin{figure}[htbp]
    \centering
    \includegraphics[width=0.8\textwidth]{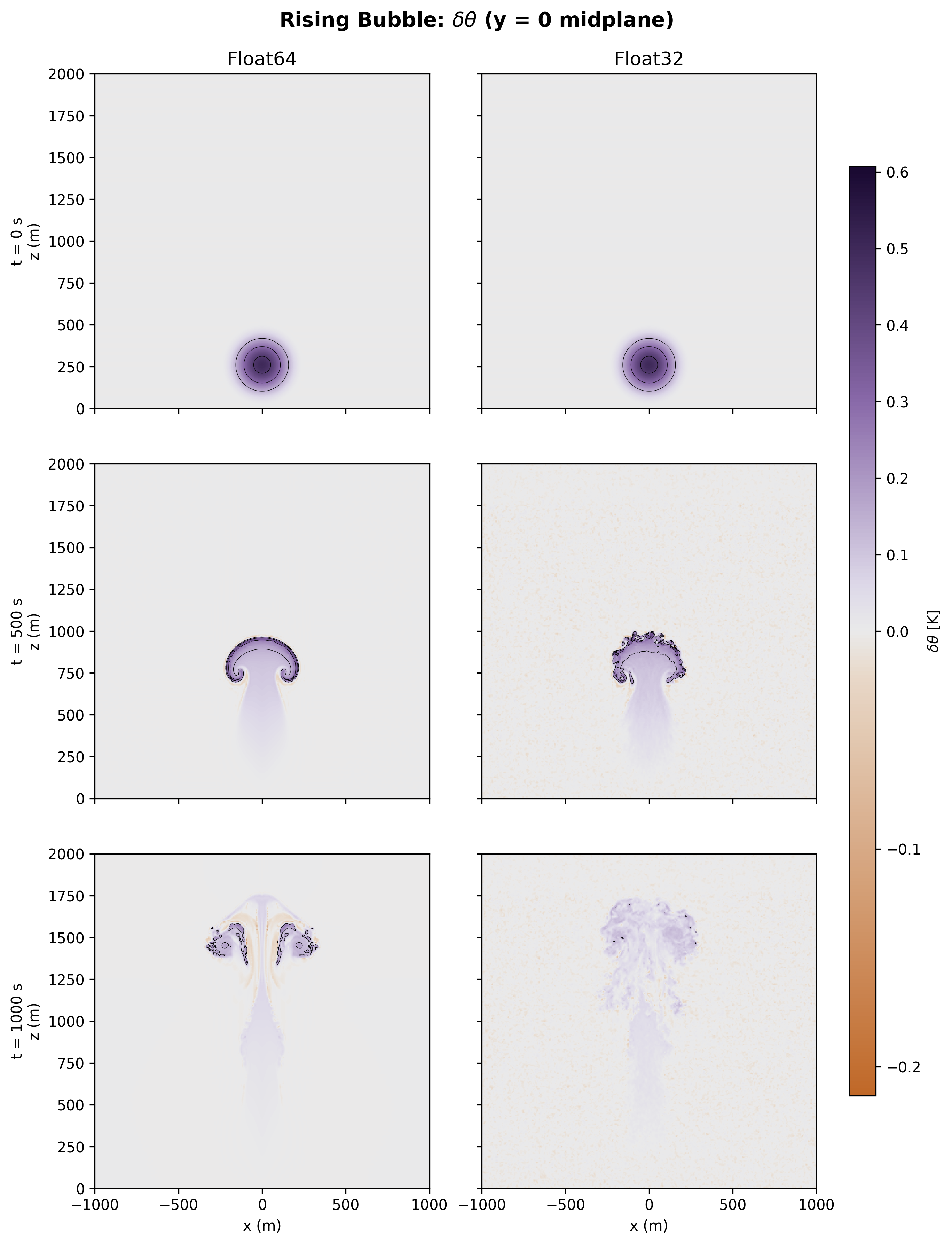}
    \caption{Slices of simulation results for potential temperature perturbation at $y=0$ with polynomial order $N=4$, refinement level $L=6$, and one coarse element in 64-bit precision (left) and 32-bit precision (right). This simulation uses the smooth initial perturbation defined in (\ref{eq:smooth-rb-ic}) for the rising thermal bubble.}
    \label{fig:rising-bubble-floatcomp-smooth}
\end{figure}
\begin{table}[htbp]
    \centering
    \caption{Timing results for the rising bubble simulation for polynomial order $N=4$ and refinement level $L=6$ (32,768,000 DOF) using eight GPUs. 
    }
    \label{tab:rb-timing}
    \begin{tabular}{llcc}
        \toprule
        Perturbation & Precision & Total Time [s] & Speedup \\
        \midrule
        \multirow{2}{*}{Sharp (\ref{eq:sharp-rb-ic})}
            & 64-bit & 6566.95 & --- \\
            & 32-bit & 3231.31 & $2.03\times$ \\
        \midrule
        \multirow{2}{*}{Smooth (\ref{eq:smooth-rb-ic})}
            & 64-bit & 6555.38 & --- \\
            & 32-bit & 3233.32 & $2.03\times$ \\
        \bottomrule
    \end{tabular}
\end{table}

\subsection{Baroclinic Instability}
\begin{figure}[htbp]
    \centering
     \begin{subfigure}[t]{0.5\textwidth}
         \centering
    \includegraphics[width=\textwidth]{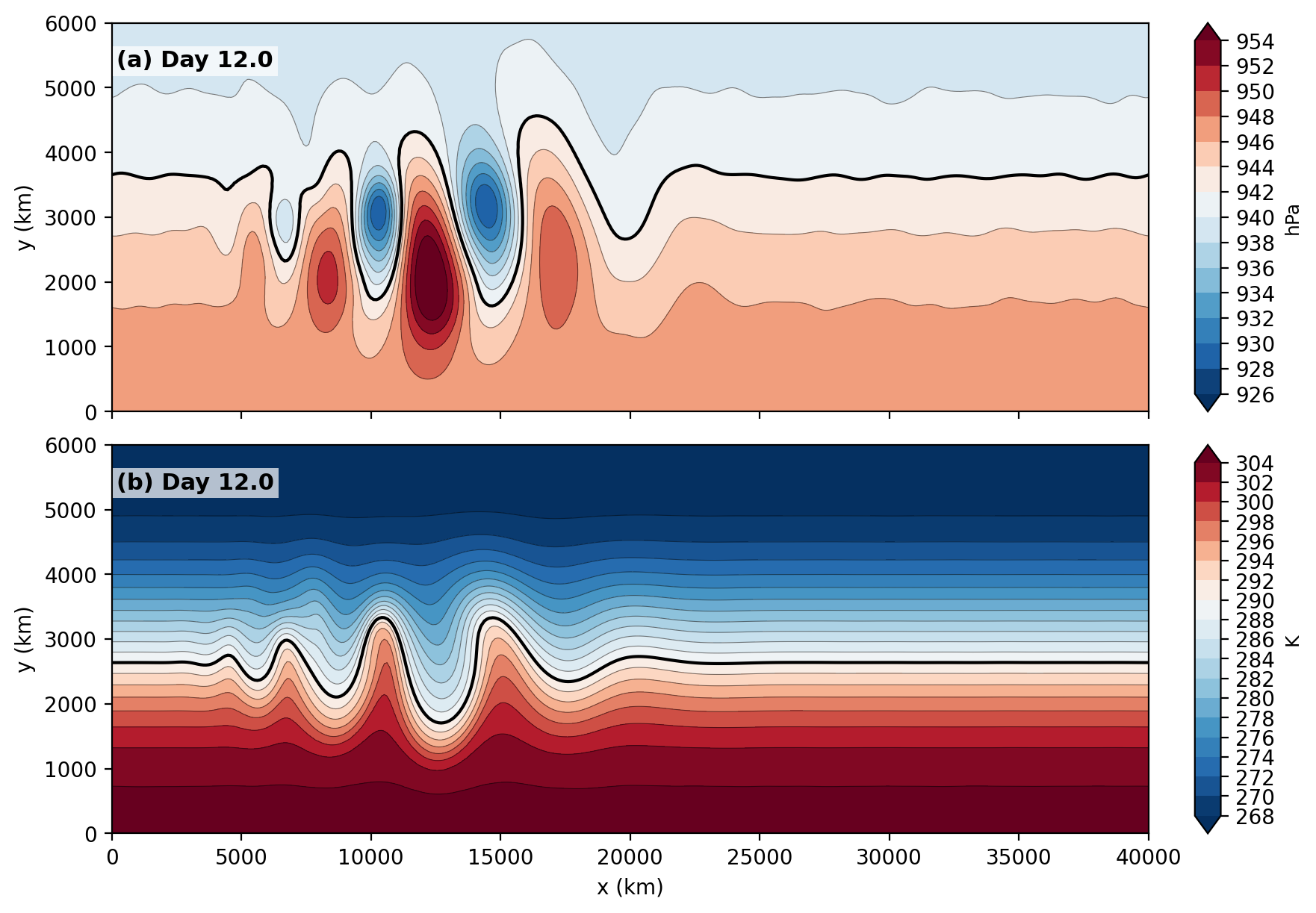}         
    \caption{ESDG}
    \label{fig:baroclinic-pres-day12-ESDG}
    \end{subfigure}%
    \hfill
    \begin{subfigure}[t]{0.5\textwidth}
         \centering
    \includegraphics[width=\textwidth]{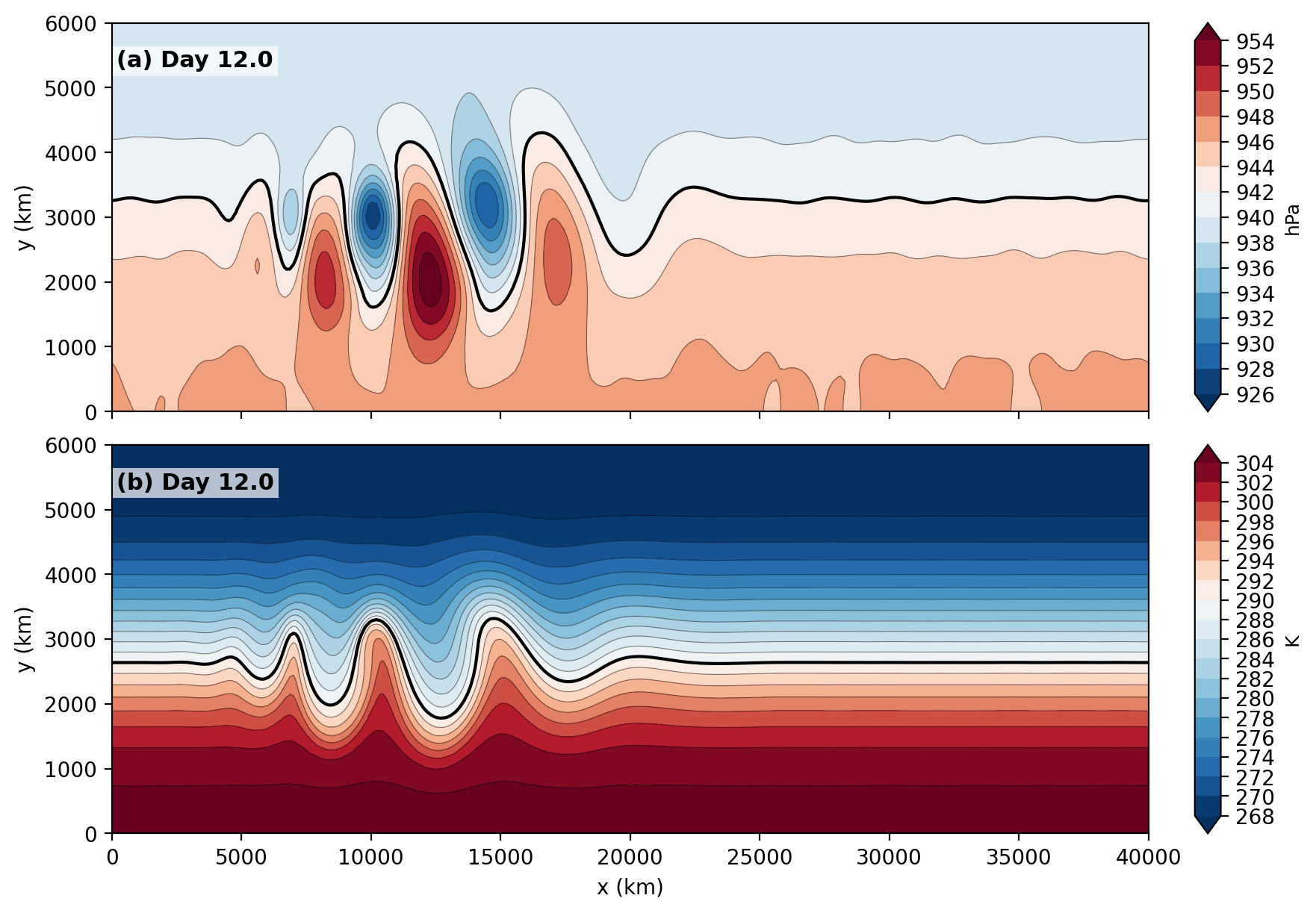} 
    \caption{NUMA}
    \label{fig:baroclinic-pres-day12-NUMA}
    \end{subfigure}%

    \caption{Slice at $z = 500~\text{m}$ of pressure (top) and temperature (bottom) on day 12 for (a) ESDG and (b) NUMA for the baroclinic instability. The simulation was run with polynomial order $N=4$, refinement level $L =3$, and $K_x = 12$, $K_y = 2$, $K_z = 1$ coarse elements in their respective directions. 
    }
    \label{fig:baroclinic-pres-day12-combined}
\end{figure}

We run a 3D baroclinic wave test case following the benchmark from \cite{Ullrich2015}.
The domain is $x \in [0, 40000]~\text{km}$ (periodic), $y \in [0, 6000]~\text{km}$ (reflecting), and $z \in [0, 30]~\text{km}$ (reflecting).
The background state is a balanced zonal wind, 
to which a Gaussian perturbation is then added, triggering an instability.
We run this simulation
with polynomial order \(N=4\), refinement level \(L=3\),
and $K_x = 12$, $K_y = 2$, and $K_z = 1$ coarse elements in their respective directions.
\RevThree{This results in 1,536,000 DOF and 
horizontal resolutions of approximately 83 km in \(x\) and 75 km in \(y\), and a vertical resolution of 750 m in \(z\).}
We run this with a time step of $\Delta t = 2.0~\text{s}$ (Courant number $C \approx 1.05$). 
We plot these results in Fig.~\ref{fig:baroclinic-pres-day12-combined} and see that they are qualitatively the same as those from \cite{Ullrich2015}. 
We also compare these results with the NUMA model \cite{Giraldo2013, Giraldo2024}. The NUMA result shown in Fig.~\ref{fig:baroclinic-pres-day12-NUMA} uses a different equation set (potential temperature for thermodynamics), spectral elements, \RevTwo{horizontally explicit vertically implicit} (HEVI) 
time integration using a second-order additive Runge–Kutta method, and fourth-order hyper-viscosity; in sum, it represents the best solution we can provide within our collection of existing and published codes. The ESDG results compare well with the NUMA results even though the ESDG method neither uses artificial dissipation nor benefits from a highly stable time-stepping approach.

\begin{figure}[htbp]
    \centering
    \includegraphics[width=1.0\textwidth]{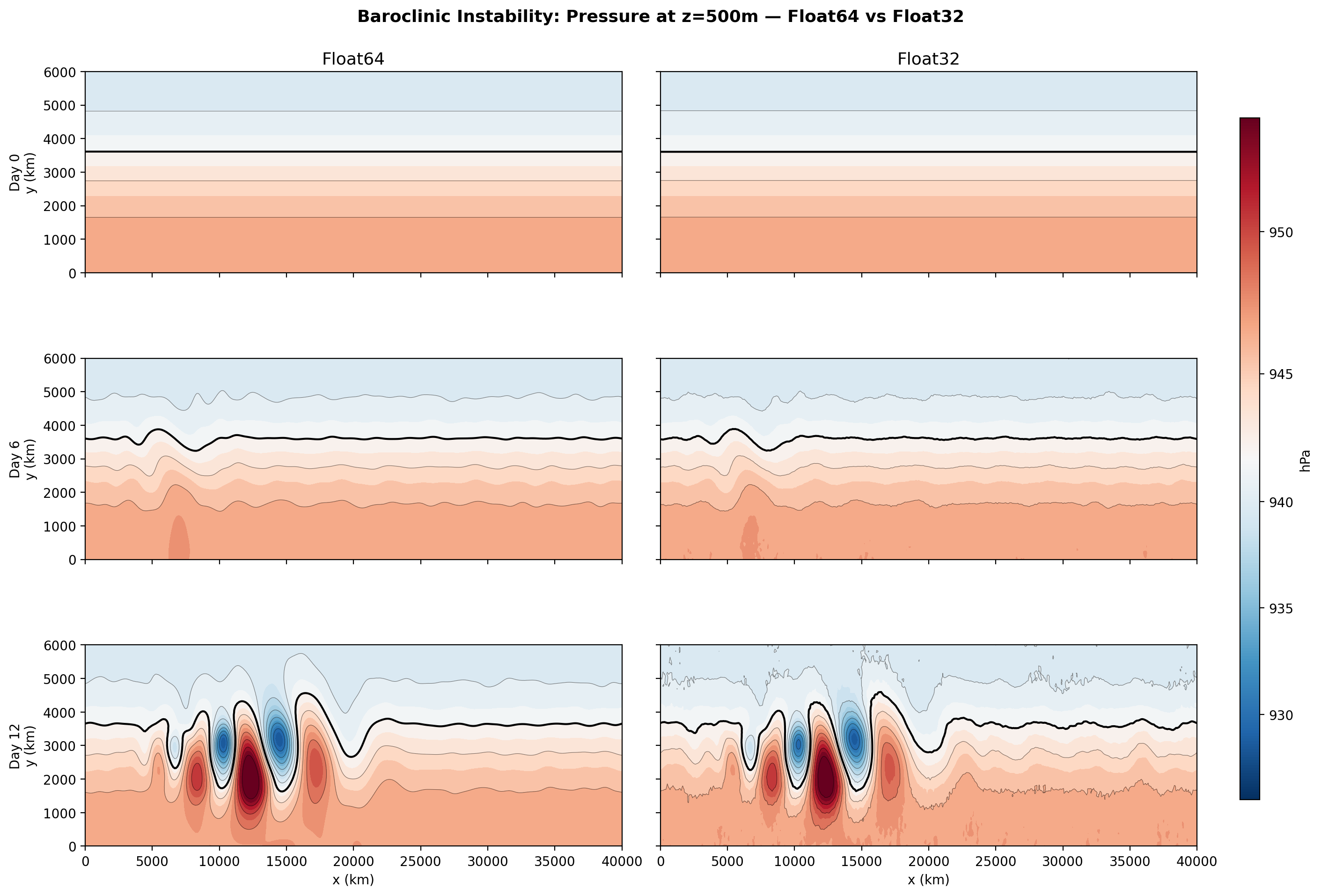}
    \caption{
    Horizontal slices of pressure at $z=500$ m at days
    0, 6, and 12 in 64-bit precision (left) and 32-bit precision (right) for the baroclinic instability. The simulation was run with polynomial order $N=4$, refinement level $L=3$, and $K_x = 12$, $K_y = 2$, $K_z = 1$ coarse elements in their respective directions.}
    \label{fig:baroclinic-pres-float-compare}
\end{figure}

\begin{figure}[htbp]
    \centering
    \includegraphics[width=1.0\textwidth]{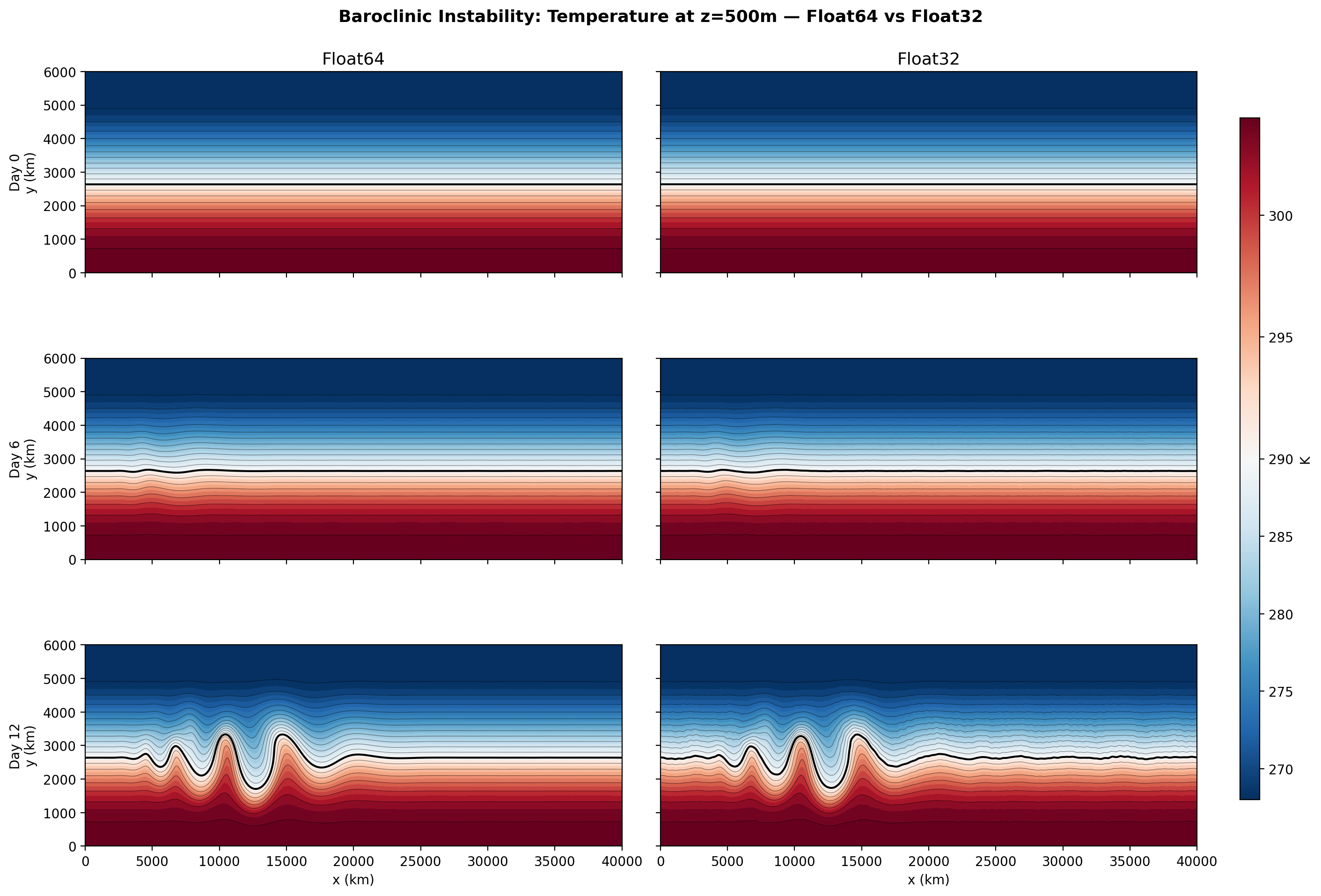}
    \caption{
    Horizontal slices of temperature at $z = 500$ m at days 0, 6, and 12 in 64-bit precision (left) and 32-bit precision (right) for the baroclinic instability. The simulation was run with polynomial order $N=4$, refinement level $L=3$, and $K_x = 12, K_y = 2, K_z = 1$ coarse elements in their respective directions.}
    \label{fig:baroclinic-temp-float-compare}
\end{figure}

As in the rising bubble case, we compare performance in 32-bit and 64-bit precision.
We summarize the timing results in Table~\ref{tab:baro-timing} and the visual results in Figs.~\ref{fig:baroclinic-pres-float-compare} and~\ref{fig:baroclinic-temp-float-compare}.
We see that the 32-bit precision solution captures the same dynamics as the 64-bit precision solution with a small amount of noise, which we attribute to the large magnitude of $\rho e$.
Global conservation of mass and total energy is verified by tracking $\sum_{K} \int_{K} \rho \, dV$ and $\sum_{K} \int_{K} \rho e \, dV$ over time; relative variations remain small for 32-bit precision  ($\sim 10^{-6}$) and 64-bit precision ($\sim 10^{-15}$).
\RevThree{
The relative change in integrated entropy from the initial to the final time step,
$(\eta_f-\eta_0)/|\eta_0|$, is approximately $-1.01\times10^{-6}$ in 32-bit precision and $-1.05\times10^{-10}$ in 64-bit precision; since the method is entropy dissipative (not entropy conserving), we do not expect to conserve entropy.}
Despite these well-behaved global diagnostics, we see visual noise in the 32-bit pressure plot.
This is most likely due to 
the use of full state variables rather than perturbation variables.
Carrying perturbation variables 
would afford the 32-bit computation greater accuracy
in small relative changes in the conserved variables, but the trade-off is additional storage and arithmetic for handling the reference state.

Table~\ref{tab:baro-timing}
shows that there is a tradeoff between speed and accuracy---we get a $1.29 \times$ speedup going from 64-bit precision to 32-bit precision.
This is less than the $2.03 \times$ speedup we saw in the rising bubble case, primarily because of the ratio of computation to communication.
Both simulations run with eight GPUs, but the rising bubble case has a DOF count of 32,768,000, whereas the baroclinic instability has a significantly smaller DOF count of 1,536,000. 
This means that the communication overhead is a larger fraction of the total runtime, and therefore 
computational optimizations such as reduced precision will have less of an overall effect.
Increasing the number of elements or the polynomial order could improve GPU utilization and further increase the speedup.
We confirm this by increasing the size of the problem to $N=4$, $L=3$, $K_x = 36$, $K_y=6$, $K_z = 3$ (41,472,000 DOFs), decreasing the time step to $\Delta t = 0.5~\text{s}$, and running for only one simulation day.
The coarse and fine timing results are listed in Table~\ref{tab:baro-timing}.
The fine baroclinic case achieves a $1.99 \times$ speedup, close to the 2.03× speedup in the rising-bubble cases.



\begin{table}[htbp]
    \centering
    \caption{Timing results for the baroclinic instability problem with both coarse (for 12~days) and fine (for a single day) configurations run on eight GPUs.
    }
    \label{tab:baro-timing}
    \begin{tabular}{llcc}
        \toprule
        Size & Precision & Total Time [s] & Speedup \\
        \midrule
        \multirow{2}{*}{Coarse}
            & 64-bit & 1886.53 & --- \\
            & 32-bit & 1456.88 & $1.29\times$ \\
        \midrule
        \multirow{2}{*}{Fine}
            & 64-bit & 11834.85 & --- \\
            & 32-bit & 5952.72 & $1.99\times$ \\
        \bottomrule
    \end{tabular}
\end{table}

%% file: sections/Conclusions.tex
\section{Conclusions} 
\label{sec:conclusions}

We presented a GPU implementation of the entropy-stable discontinuous Galerkin method for compressible Euler equations with buoyancy. We demonstrated the performance of the solver on two three-dimensional problems, showing that on NVIDIA A100 hardware, the solver achieves nearly 70\% of 64-bit peak performance for the most computationally expensive kernel (volume terms) and that the optimized volume kernel significantly reduces the computational overhead typically incurred by two-point entropy-stable fluxes.
We also demonstrated strong and weak scaling
of the solver and 
compared it with a highly optimized CPU code, showing that the GPU
kernels are $10\times$ faster and $13\times$ more energy efficient than the CPU code. We discussed different modifications implemented to reach this level of performance and reported on the performance gains of each implementation 
strategy, including reductions in complex operations and memory traffic and improvements in load balancing.
We only used explicit time integration because we wanted to focus on the performance of the ESDG spatial discretization; however, for a practical atmospheric model it is imperative, at a minimum, to 
extend the solver by incorporating horizontally explicit vertically implicit time-integration methods (see, e.g., \cite{Souza2023, Giraldo2024}).